\documentclass{amsart}
\usepackage{amssymb, latexsym}

\newcommand\card{\operatorname{card}}

\newcommand \Ivl {\left( [0, \, \tau] \right) \, }

\theoremstyle{plain}
\newtheorem{thm}{Theorem}

\newtheorem{prop}[thm]{Proposition}
\newtheorem{lem}[thm]{Lemma}

\numberwithin{equation}{section} \numberwithin{thm}{section}

\begin{document}

\title[Radial defocusing cubic wave equation ]{Global well-posedness for the radial defocusing cubic wave equation
on $\mathbb{R}^{3}$ and for rough data }

\author{Tristan Roy}
\address{Nagoya University}
\email{tristanroy@math.nagoya-u.ac.jp}

\vspace{-0.3in}

\begin{abstract}
We prove global well-posedness for the radial defocusing cubic wave
equation

\begin{equation}
\left\{
\begin{array}{ccl}
\partial_{tt} u  - \Delta u & = & -u^{3} \\
u(0,x)& = & u_{0}(x)    \\
\partial_{t} u(0,x) & = & u_{1}(x)
\end{array}
\right. \nonumber
\end{equation}
with data $\left( u_{0}, \, u_{1} \right) \in H^{s} \times H^{s-1}$,
$1 > s >\frac{7}{10}$. The proof relies upon a Morawetz-Strauss-type
inequality that allows us to control the growth of an almost
conserved quantity.

\end{abstract}

\maketitle

\section{Introduction}

We shall study the defocusing cubic wave equation on
$\mathbb{R}^{3}$

\begin{equation}
\left\{
\begin{array}{ccl}
\partial_{tt} u  - \Delta u & = & -u^{3} \\
u(0,x)& = & u_{0}(x)    \\
\partial_{t} u(0,x) & = & u_{1}(x)
\end{array}
\right. \label{Eqn:WaveEqRad}
\end{equation}
We shall focus on the strong solutions of the defocusing cubic wave
equation on some interval $[0,T]$ i.e real-valued maps
$(u,\partial_{t} u) \in C \left([0, \, T], \, H^{s}(\mathbb{R}^{3})
\right) \times C \left( [0, \,T], \, H^{s-1} ( \mathbb{R}^{3})
\right)$ that satisfy for $t \in [0, \, T] $ the following integral
equation

\begin{equation}
\begin{array}{ll}
u(t) & = \cos(tD) u_{0} + D^{-1} \sin(tD) u_{1} - \int_{0}^{t}
D^{-1} \sin \left( (t-t^{'}) D \right) u^{3}(t^{'}) \, dt^{'}
\end{array}
\end{equation}
with $(u_{0},u_{1})$ lying in $H^{s} \times H^{s-1}$. Here $H^{s}$
is the usual inhomogeneous Sobolev space i.e $H^{s}$ is the
completion of the Schwartz space $\mathcal{S}(\mathbb{R}^{3})$ with
respect to the norm

\begin{equation}
\begin{array}{ll}
\| f \|_{H^{s}} & :=  \| (1+ D^{s}) f \|_{L^{2}(\mathbb{R}^{3})}
\end{array}
\end{equation}
where $D$ is the operator defined by

\begin{equation}
\begin{array}{ll}
\widehat{Df}(\xi) & := |\xi| \hat{f}(\xi)
\end{array}
\end{equation}
and $\hat{f}$ denotes the Fourier transform
\begin{equation}
\begin{array}{ll}
\hat{f}(\xi) & := \int_{\mathbb{R}^{3}} f(x) e^{-i x \cdot \xi} \,
dx
\end{array}
\end{equation}
Here $H^{s} \times H^{s-1}$ is the product space of $H^{s}$ and
$H^{s-1}$ endowed with the standard norm $ \| (f,g) \|_{H^{s} \times
H^{s-1}} := \| f \|_{H^{s}} + \| g \|_{H^{s-1}}$.

It is known \cite{lindsogge} that (\ref{Eqn:WaveEqRad}) is locally
well-posed in $H^{s}(\mathbb{R}^{3}) \times H^{s-1}(\mathbb{R}^{3})$
for $s \geq \frac{1}{2}$. Moreover if $s > \frac{1}{2}$ the time of
local existence only depends on the norm of the initial data $\|
(u_{0},u_{1}) \|_{H^{s} \times H^{s-1}}$.

Now we turn our attention to the global well-posedness theory of
(\ref{Eqn:WaveEqRad}). In view of the above local well-posedness
theory and standard limiting arguments it suffices to establish an a
priori bound of the form

\begin{equation}
\begin{array}{ll}
\| u(T) \|_{H^{s}} + \| \partial_{t} u (T) \|_{H^{s-1}} & \leq  C
\left( s, ( \| u_{0}\|, \| u_{1} \| )_{H^{s} \times H^{s-1}}, T
\right)
\end{array}
\end{equation}
for all times $0 < T < \infty $ and all smooth-in-time
Schwartz-in-space solutions $(u,\partial_{t} u): [0, \, T] \times
\mathbb{R}^{3} \rightarrow \mathbb{R}$, where the right-hand side is
a finite quantity depending only on $s$, $ \| u_{0} \|_{H^{s}} $, $
\| u_{1} \|_{H^{s-1}}$ and $T$. Therefore in the sequel we shall
restrict attention to such smooth solutions.

The defocusing cubic wave equation (\ref{Eqn:WaveEqRad}) enjoys the
following energy conservation law

\begin{equation}
\begin{array}{ll}
E(u(t)) & := \frac{1}{2} \int_{\mathbb{R}^{3}} (\partial_{t}
u)^{2}(x,t) \, dx  + \frac{1}{2} \int_{\mathbb{R}^{3}} | D u
(x,t)|^{2} \, dx + \frac{1}{4} \int_{\mathbb{R}^{3}} u^{4}(x,t) \,
dx
\end{array}
\end{equation}
Combining this conservation law to the local well-posedness theory
we immediately have global well-posedness for (\ref{Eqn:WaveEqRad})
and for $s=1$.

In this paper we are interested in studying global well-posedness
for (\ref{Eqn:WaveEqRad}) and for data below the energy norm, i.e
$s<1$. It is conjectured that (\ref{Eqn:WaveEqRad}) is globally
well-posed in $H^{s} (\mathbb{R}^{3}) \times H^{s-1}(\mathbb{R}^{3})
$ for all $s > \frac{1}{2}$. The global existence for the defocusing
cubic wave equation has been the subject of several papers. Let us
some mention some results for data lying in a slightly different
space than $H^{s} \times H^{s-1}$ i.e $\dot{H}^{s} \times
\dot{H}^{s-1}$ \footnote{More precisely the data lie in $\dot{H}^{s} \cap L^{4} \times \dot{H}^{s-1}$}. Here $\dot{H}^{s}$ is the usual homogeneous Sobolev
space i.e the completion of Schwartz functions $\mathcal{S} \left(
\mathbb{R}^{3} \right)$ with respect to the norm

\begin{equation}
\begin{array}{ll}
\| f \|_{\dot{H}^{s}} & = \| D^{s} f \|_{L^{2} \left( \mathbb{R}^{3}
\right)}
\end{array}
\end{equation}

Kenig, Ponce and Vega \cite{kenponcevega} were the first to prove
that (\ref{Eqn:WaveEqRad}) is globally well-posed for $1 >s
> \frac{3}{4}$. They used the \emph{Fourier truncation method} discovered by
Bourgain \cite{bourg}. I. Gallagher and F. Planchon
\cite{gallagplanch}  proposed a different method to prove global
well-posedness for $1 >s
> \frac{3}{4}$. H. Bahouri and Jean-Yves
Chemin \cite{bahchemin} proved global-wellposedness for
(\ref{Eqn:WaveEqRad}) and for $s=\frac{3}{4}$ by using a non linear
interpolation method and logarithmic estimates from S. Klainermann
and D. Tataru \cite{klaintat}. We shall consider global
well-posedness for the radial defocusing cubic wave equation i.e
global existence for the initial value problem (\ref{Eqn:WaveEqRad})
with radial data. The main result of this paper is the following one

\begin{thm}
The radial defocusing cubic wave equation is globally well-posed in
$ H^{s} \times H^{s-1} $ for $ 1 > s > \frac{7}{10} $. Moreover if
$T$ large then

\begin{equation}
\begin{array}{ll}
\| u(T) \|^{2}_{H^{s}} + \| \partial_{t} u(T) \|^{2}_{H^{s-1}} &
\leq C \left( \| u_{0} \|_{H^{s}}, \, \| u_{1} \|_{H^{s-1}} \right)
T^{\frac{16s-10}{10s-7}+}
\end{array}
\label{Eqn:LgEstRaduT1}
\end{equation}
for $ \frac{5}{6} \geq s >  \frac{7}{10}$ and

\begin{equation}
\begin{array}{ll}
\| u(T) \|^{2}_{H^{s}} + \| \partial_{t} u(T) \|^{2}_{H^{s-1}} &
\leq C \left( \| u_{0} \|_{H^{s}}, \, \| u_{1} \|_{H^{s-1}} \right)
T^{\frac{2s}{2s-1}+}
\end{array}
\label{Eqn:LgEstRaduT2}
\end{equation}
for $ 1 > s > \frac{5}{6}$. Here $C \left( \| u_{0} \|_{H^{s}}, \,
\| u_{1} \|_{H^{s-1}} \right)$ is a constant only depending on $ \|
u_{0} \|_{H^{s}} $ and $\| u_{1} \|_{H^{s-1}}$.
\label{Thm:GwpRad710}
\end{thm}

We set some notation that appear throughout the paper. Given $A,B$
positive number $A \lesssim B$ means that there exists a universal
constant $K$ such that $A \leq K B$. We say that $K_{0}$ is the
constant determined by the relation $A \lesssim B$ if $K_{0}$ is the
smallest $K$ such that $A \leq K B$ is true. We write $A \sim B$
when $A \lesssim B$ and $B \lesssim A$. $A << B$ denotes $ A \leq K
B$ for some universal constant $K < \frac{1}{100}$ . We also use the
notations $A+ = A + \epsilon$, $A-=A - \epsilon$ for some universal
constant $0 < \epsilon << 1$. Let $\nabla$ denote the gradient
operator. If $J$ is an interval then $|J|$ is its size. If $E$ is a
set then $\card(E)$ is its cardinal. Let $I$ be the following
multiplier

\begin{equation}
\begin{array}{ll}
\widehat{If}(\xi) & := m(\xi) \hat{f}(\xi)
\end{array}
\end{equation}
where $m(\xi): =  \eta \left( \frac{\xi}{N} \right)$, $\eta$ is a
smooth, radial, nonincreasing in $|\xi|$ such that

\begin{equation}
\begin{array}{ll}
\eta (\xi) & :=  \left\{
\begin{array}{l}
1, \, |\xi| \leq 1 \\
\left( \frac{1}{|\xi|} \right)^{1-s}, \, |\xi| \geq 2
\end{array}
\right.
\end{array}
\end{equation}
and $N>>1$ is a dyadic number playing the role of a parameter to be
chosen. We shall abuse the notation and write $m (|\xi|)$ for
$m(\xi)$, thus for instance $m(N)=1$.

We recall some basic results regarding the defocusing cubic wave
equation. Let $\lambda \in \mathbb{R}$ and $u_{\lambda}$ denote the
following function

\begin{equation}
\begin{array}{ll}
u_{\lambda}(t,x) & := \frac{1}{\lambda} u \left( \frac{t}{\lambda},
\frac{x}{\lambda} \right)
\end{array}
\end{equation}
If $u$ satisfies (\ref{Eqn:WaveEqRad}) with data $(u_{0},u_{1})$
then $u_{\lambda}$ also satisfies (\ref{Eqn:WaveEqRad}) but with
data $ \left( \frac{1}{\lambda} u_{0} \left( \frac{x}{\lambda}
\right), \frac{1}{\lambda^{2}} u_{1} \left( \frac{x}{\lambda}
\right) \right)$. If $u$ satisfies the radial defocusing cubic wave
equation then $u$ is radial.

Now we recall some standard estimates that we use later in this
paper.

\begin{prop} {\textbf{''Strichartz estimates in 3 dimensions''} (See \cite{ginebvelo},
\cite{lindsogge})}. Let $m \in [0, \,1]$. If $u$ is a strong
solution to the IVP problem

\begin{equation}
\left\{
\begin{array}{l}
 \partial_{tt} u - \Delta   u = F \\
 u(0,x) = f(x) \in \dot{H}^{m} \\
 \partial_{t} u(0,x) =  g(x) \in \dot{H}^{m-1} \\
\end{array}
\right. \label{Eqn:LinearWave}
\end{equation}
then we have for $0 \leq \tau < \infty$

\begin{equation}
\begin{array}{ll}
\| u \|_{L_{t}^{q}\left( [0, \, \tau] \right) L_{x}^{r}} + \| u
\|_{C \left( [0, \, \tau]; \dot{H}^{m} \right)} + \| \partial_{t} u
\|_{C \left( [0, \, \tau]; \dot{H}^{m-1} \right)} & \lesssim \| f
\|_{\dot{H}^{m}} + \| g \|_{\dot{H}^{m-1}} + \| F
\|_{L_{t}^{\tilde{q}} \left( [0, \, \tau] \right) L_{x}^{\tilde{r}}}
\end{array}
\nonumber
\end{equation}
under two assumptions

\begin{itemize}

\item $(q,r)$ lie in the set $\mathcal{W}$ of \textit{wave-admissible} points i.e

\begin{equation}
\begin{array}{ll}
\mathcal{W}  & := \left\{ (q, \, r): ( q, \, r ) \in (2, \,\infty ]
\times [2,\infty), \, \frac{1}{q}+\frac{1}{r} \leq \frac{1}{2}
\right\}
\end{array}
\label{Eqn:StrCondition1}
\end{equation}

\item $(\tilde{q}, \tilde{r})$ lie in the dual set $\mathcal{W}^{'}$ of $\mathcal{W}$ i.e

\begin{equation}
\begin{array}{ll}
\mathcal{W}^{'}  & :=  \left\{ (\tilde{q}, \tilde{r}):
\frac{1}{\tilde{q}} + \frac{1}{q}=1, \, \frac{1}{r}+
\frac{1}{\tilde{r}}=1, \, (q,r) \in \mathcal{W} \right\}
\end{array}
\end{equation}

\item $(q,r,\tilde{q}, \tilde{r})$ satisfy the \textit{dimensional analysis}
conditions

\begin{equation}
\begin{array}{l}
\frac{1}{q} + \frac{3}{r} =  \frac{3}{2} -m
\end{array}
\label{Eqn:StrCondition2}
\end{equation}
and

\begin{equation}
\begin{array}{l}
\frac{1}{\tilde{q}} + \frac{3}{\tilde{r}} -2 = \frac{3}{2} -s
\end{array}
\label{Eqn:StrCondition3}
\end{equation}

\end{itemize}

\end{prop}
We also have the well-known estimate

\begin{prop}{\textbf{''Radial Sobolev inequality''}}
If $u:\mathbb{R}^{3} \rightarrow \mathbb{C}$ is radial and smooth
\begin{eqnarray}
|u(x)| & \lesssim & \frac{ \|u\|_{\dot{H}^{1}}}{|x|^{\frac{1}{2}}}
\label{Eqn:SobolevRad}
\end{eqnarray}
\end{prop}
The Hardy-type inequality is proved in \cite{caz}

\begin{prop}{\textbf{''Hardy-type inequality''}}
If $1<p<3$  and $u:\mathbb{R}^{3} \rightarrow \mathbb{C}$ is smooth

\begin{equation}
\begin{array}{ll}
\| \frac{u}{|x|} \|_{L^{p}} & \leq  \frac{3}{3-p} \| Df \|_{L^{p}}
\end{array}
\label{Eqn:HardyIneq}
\end{equation}
\end{prop}

Some variables appear frequently in this paper. We define them now.

We say that $(q,r)$ is a $m$-wave admissible pair if $0 \leq  m \leq
1$ and $(q,r)$ satisfy the two following conditions

\begin{itemize}
\item $(q,r) \in \mathcal{W}$

\item $ \frac{1}{q} + \frac{3}{r} =  \frac{3}{2}-m $

\end{itemize}

Let $J=[a, \,b]$ be an interval included in $[0, \, \infty)$. Let
$Z_{m,s}(J)$ denote the following number

\begin{equation}
\begin{array}{ll}
Z_{m,s}(J) & :=  \sup_{q,r} \left( \| D^{1-m} Iu
\|_{L_{t}^{q}(J)L_{x}^{r}} + \| D^{-m} I \partial_{t} u
\|_{L_{t}^{q}(J)L_{x}^{r}} \right)
\end{array}
\end{equation}
where the $\sup$ is taken over $m$-wave admissible $(q,r)$  and let

\begin{equation}
\begin{array}{ll}
Z(J) & :=  \sup_{m \in [0, \, 1)} Z_{m,s}(J)
\end{array}
\end{equation}
Let

\begin{equation}
\begin{array}{ll}
R_{1}(J) & :=  \int_{J} \int_{\mathbb{R}^{3}} \frac{\nabla I
u(t,x).x}{|x|} \left( (Iu)^{3}(t,x) -Iu^{3}(t,x) \right) \,  dx dt
\end{array}
\end{equation}
and

\begin{equation}
\begin{array}{ll}
R_{2}(J) & := \int_{J} \int_{\mathbb{R}^{3}} \frac{Iu(t,x)}{|x|}
\left( (Iu)^{3}(t,x) -Iu^{3}(t,x) \right) \,  dx dt
\end{array}
\end{equation}
If $J=[0, \tau]$ we shall abuse the notation and write

\begin{equation}
\begin{array}{ll}
Z(\tau) & := Z \left( J \right) \\
R(\tau) & := R \left( J \right)
\end{array}
\nonumber
\end{equation}

Some estimates that we establish throughout the paper require a
Paley-Littlewood decomposition. We set it up now. Let $\phi(\xi)$ be
a real, radial, nonincreasing function that is equal to $1$ on the
unit ball $\left\{ \xi \in \mathbb{R}^{3}: \, |\xi| \leq 1 \right\}$
and that that is supported on $\left\{ \xi \in \mathbb{R}^{3}: \,
|\xi| \leq 2 \right\}$. Let $\psi$ denote the function

\begin{equation}
\begin{array}{ll}
\psi(\xi) & := \phi(\xi) - \phi(2 \xi)
\end{array}
\end{equation}
If $M \in 2^{\mathbb{Z}}$ is a dyadic number we define the
Paley-Littlewood operators in the Fourier domain by

\begin{equation}
\begin{array}{ll}
\widehat{P_{\leq M} f}(\xi) & := \phi \left( \frac{\xi}{M} \right)
\hat{f}(\xi) \\
\widehat{P_{M} f}(\xi) & := \psi \left( \frac{\xi}{M} \right)
\hat{f}(\xi) \\
\widehat{P_{> M} f}(\xi) & := \hat{f}(\xi) - \widehat{P_{\leq M}
f}(\xi)
\end{array}
\end{equation}
Since $\sum_{M \in  2^{\mathbb{Z}}} \psi \left( \frac{\xi}{M}
\right)=1$ we have

\begin{equation}
\begin{array}{ll}
f & = \sum_{M \in 2^{\mathbb{Z}}} P_{M} f
\end{array}
\end{equation}

We conclude this introduction by giving the main ideas of the proof
of theorem \ref{Thm:GwpRad710} and explaining how the paper is
organized. Following the proof of the global well-posedness for
$s=1$ we try to compare for every $T>0$ the relevant quantity  $ \|
\left( u(T),
\partial_{t} u(T) \right) \|_{H^{s} \times H^{s-1}}$ to the supremum
of the energy conservation law $ \sup_{t \in [0, \, T]} E \left(
u(t) \right)$. Unfortunately this strategy does not work if $s<1$
since the energy can be infinite. We get around this difficulty by
using the $I$-method designed by J. Colliander, M. Keel, G.
Staffilani, H.Takaoka and T. Tao \cite{almckstt} and successfully
applied to prove global well-posedness for semilinear Schr\"odinger
equations and for rough data. The idea consists of introducing the
following smoothed energy

\begin{equation}
\begin{array}{ll}
E \left( Iu(t) \right) & := \frac{1}{2} \int_{\mathbb{R}^{3}} \left|
\partial_{t} I u(x,t) \right|^{2} \, dx + \frac{1}{2} \int_{\mathbb{R}^{3}} | D I u(x,t)|^{2} \, dx +
\frac{1}{4} \int_{\mathbb{R}^{3}} |I u(x,t)|^{4} \, dx
\end{array}
\end{equation}
We prove in section \ref{sec:RelHsNrj} that $\| \left( u(T), \,
\partial_{t} u(T) \right) \|^{2}_{H^{s} \times H^{s-1}}$ and the
supremum of the smoothed energy on $[0, \, T]$ are comparable.
Therefore we try to estimate $\sup_{t \in [0, \, T]}E \left ( Iu(t)
\right)$ in order to give an upper bound of $ \| ( u(T), \,
\partial_{t} u(T)) \|_{H^{s} \times H^{s-1}} $. For convenience we place the
mollified energy at time zero into $[0, \, \frac{1}{2}]$ by choosing
the right scaling factor $\lambda$. This operation shows that we are
reduced to estimate $\sup_{t \in [0, \, \lambda T]} E \left( I
u_{\lambda}(t) \right)$. In section \ref{sec:LocalBd} we prove that
we can locally control a variable namely $Z(J)$ provided that the
interval $J$ satisfies some constraints that give some information
about its size. $\sup_{t \in J} E \left( I u_{\lambda}(t) \right) $
is estimated by the fundamental theorem of calculus. The upper bound
depends on the parameter $N$ and the controlled quantity $Z(J)$.
This estimate is established in section \ref{sec:AlmCon}. Now we can
iterate: the process generates a sequence of intervals $\left( J_{i}
\right)$ that cover the whole interval $[0, \, \lambda T]$ and
satisfy the same constraints as $J$. We should be able to estimate
$\sup_{t \in [0, \, \lambda T]}E \left ( Iu_{\lambda}(t) \right)$
provided that we can control the number of intervals $J_{i}$. This
requires the establishment of a long time estimate, the so-called
almost Morawetz-Strauss inequality. This estimate is proved in
section \ref{sec:AlmMorIneq}. It depends on some remainder integrals
that are estimated in section \ref{subsec:RemainderMor}. Combining
this inequality to the radial Sobolev inequality
(\ref{Eqn:SobolevRad}) we can give an upper bound of the cardinal of
$\left( J_{i} \right)$. The proof of theorem \ref{Thm:GwpRad710} is
given in section \ref{sec:PfGwRad710}.

\vspace{5mm}

$\textbf{Acknowledgements}:$ The author would like to thank his
advisor Terence Tao for introducing him to this topic and is
indebted to him for many helpful conversations and encouragement
during the preparation of this paper.


\section{Proof of global well-posedness for $1 > s > \frac{7}{10}$}
\label{sec:PfGwRad710}

In this section we prove the global existence of
(\ref{Eqn:WaveEqRad}) for $1 > s > \frac{7}{10}$. Our proof relies
on some intermediate results that we prove in later sections. More
precisely we shall show the following

\begin{prop}{\textbf{'' $H^{s}$ norms and mollified energy estimates ''} }
Let $T>0$. Then

\begin{equation}
\begin{array}{ll}
\| u(T) \|^{2}_{H^{s}} + \|\partial_{t} u(T)\|^{2}_{H^{s-1}} &
\lesssim  \| u_{0} \|^{2}_{H^{s}} + \left( T^{2}+1 \right)\sup_{t
\in [0, \, T]} E \left( Iu(t) \right)
\end{array}
\label{Eqn:NrjEstim1}
\end{equation}
for every $u$. \label{Prop:NrjEst}
\end{prop}

\begin{prop}{\textbf{"Local boundedness"}}
Let $J=[a,b]$ be an interval included in $[0, \, \infty)$. Assume
that $E \left( Iu(a) \right) \leq 2$ and that $u$ satisfies
(\ref{Eqn:WaveEqRad}). There exist $C_{1}$, $C_{2}$ small and
positive constants such that if $J$ satisfies

\begin{equation}
\begin{array}{ll}
\| Iu \|_{L_{t}^{6}(J) L_{x}^{6}} & \leq
\frac{C_{1}}{|J|^{\frac{1}{3}}}
\end{array}
\label{Eqn:Condition1}
\end{equation}
and

\begin{equation}
\begin{array}{ll}
|J| & \leq  C_{2} N^{\frac{1-s}{s-\frac{1}{2}}}
\end{array}
\label{Eqn:Condition2}
\end{equation}
then  $Z(J)  \lesssim  1$.
\label{prop:LocalBdRad}
\end{prop}

\begin{prop}{\textbf{"Almost conservation law"}}
Let $J=[a,b]$ be an interval included in $[0,\infty)$. Assume that
$u$ satisfies (\ref{Eqn:WaveEqRad}). Then we have

\begin{equation}
\begin{array}{ll}
\left| \sup_{t \in J} E(Iu(t)) - E(Iu(a)) \right| & \lesssim
\frac{Z^{4} (J)}{N^{1-}} \label{Eqn:EstNrjRad}
\end{array}
\end{equation}
\label{prop:EstNrjRad}
\end{prop}

\begin{prop}{\textbf{"Almost Morawetz-Strauss inequality"}}
Let $T \geq 0$. Assume that $u$ satisfies (\ref{Eqn:WaveEqRad}).
Then we have

\begin{equation}
\begin{array}{ll}
 \int_{0}^{T} \int_{\mathbb{R}^{3}} \frac{\left| Iu \right|^{4}(t,x)}{|x|}
dx \, dt -2 \left( E \left( I u(0) \right) + E \left( Iu(T) \right)
\right)  & \lesssim  \left| R_{1}(T)\right| + \left| R_{2}(T)
\right|
\end{array}
\end{equation}
\label{prop:AlmMorawetz}
\end{prop}
and

\begin{prop}{\textbf{"Estimate of integrals"}}
Let $J$ be an interval included in $[0,\infty)$. Then if $i=1,2$ we
have

\begin{eqnarray}
R_{i}(J) & \lesssim & \frac{Z^{4}(J)}{N^{1-}}
\end{eqnarray}
\label{prop:RemainderMor}
\end{prop}

For the remainder of the section we show how proposition
\ref{prop:LocalBdRad}, \ref{prop:EstNrjRad}, \ref{prop:AlmMorawetz}
and \ref{prop:RemainderMor} imply Theorem \ref{Thm:GwpRad710}.

Let $T>0$ and $N=N(T)>>1$ be a parameter to be chosen later. There
are three steps to prove Theorem \ref{Thm:GwpRad710}.

\begin{enumerate}

\item \textbf{Scaling}. Let $\lambda >>1$ to be chosen later. Then by Plancherel theorem

\begin{equation}
\begin{array}{ll}
\| D I u_{\lambda}(0) \|^{2}_{L^{2}} & \lesssim \int_{|\xi| \leq 2N}
|\xi|^{2}  | \widehat{u_{\lambda}}(0,\xi)|^{2} \, d\xi
+ \int_{|\xi| \geq 2N} |\xi|^{2} \frac{N^{2(1-s)}}{|\xi|^{2(1-s)}} |\widehat{u_{\lambda}}(0, \xi)|^{2} \, d\xi \\
& \lesssim N^{2(1-s)} \| u_{\lambda}(0) \|^{2}_{\dot{H}^{s}} \\
& \lesssim N^{2(1-s)} \lambda^{1-2s} \| u_{0} \|^{2}_{\dot{H}^{s}} \\
& \lesssim N^{2(1-s)} \lambda^{1-2s} \| u_{0} \|^{2}_{H^{s}}
\end{array}
\label{Eqn:DIuSc}
\end{equation}

\begin{equation}
\begin{array}{ll}
\| \partial_{t} I u_{\lambda}(0) \|^{2}_{L^{2}} & \lesssim
\int_{|\xi| \leq 2N } |\widehat{\partial_{t}
u_{\lambda}}(0,\xi)|^{2} \, d\xi
+ \int_{|\xi| \geq 2N} \frac{N^{2(1-s)}}{|\xi|^{2(1-s)}} |\widehat{\partial_{t} u_{\lambda}}(0,\xi)|^{2} \, d\xi \\
& \lesssim N^{2(1-s)} \| \partial_{t} u_{\lambda}(0) \|^{2}_{H^{s-1}} \\
& \lesssim N^{2(1-s)} \left( \int_{|\xi| \leq 1}
|\widehat{\partial_{t} u_{\lambda}}(0,\xi)|^{2} \, d\xi
+ \int_{|\xi| \geq 1} |\xi|^{2(s-1)} |\widehat{\partial_{t} u_{\lambda}}(0,\xi)|^{2} \, d\xi \right) \\
& \lesssim N^{2(1-s)} \left( \frac{1}{\lambda} \int_{|\xi| \leq
\lambda} | \widehat{u_{1}}(\xi) |^{2} \, d\xi
+ \lambda^{1-2s} \int_{|\xi| \geq \lambda} |\xi|^{2(s-1)} |\widehat{u_{1}}(\xi)|^{2} \, d\xi \right) \\
& \lesssim N^{2(1-s)} \lambda^{1-2s} \| u_{1} \|^{2}_{H^{s-1}}
\end{array}
\label{Eqn:DtIuSc}
\end{equation}
By homogeneous Sobolev embedding

\begin{equation}
\begin{array}{ll}
\| I u_{\lambda}(0) \|^{2}_{L^{4}} & \lesssim \int_{\mathbb{R}^{3}}
|\xi|^{\frac{3}{2}}
|\widehat{Iu_{\lambda}}(0,\xi)|^{2} \, d\xi \\
& \lesssim \int_{|\xi| \leq 2N} |\xi|^{\frac{3}{2}}
|\widehat{u_{\lambda}}(0,\xi)|^{2} \, d\xi +
\int_{|\xi| \geq 2N} |\xi|^{\frac{3}{2}} \frac{N^{2(1-s)}}{|\xi|^{2(1-s)}} |\widehat{u_{\lambda}}(0,\xi)|^{2} \, d\xi \\
& \lesssim \frac{1}{\lambda^{\frac{1}{2}}} \int_{|\xi| \leq 2N
\lambda} |\xi|^{\frac{3}{2}} |\widehat{u_{0}}(\xi)|^{2} \, d\xi
+ N^{2(1-s)} \lambda^{\frac{3}{2}-2s}  \int_{|\xi| \geq 2N \lambda} |\xi|^{2s-\frac{1}{2}} |\widehat{u_{0}}(\xi)|^{2} \, d\xi \\
& \lesssim \frac{ \max{\left( N^{\frac{3}{2}-2s}
\lambda^{\frac{3}{2}-2s},1 \right)}}{\lambda^{\frac{1}{2}}} \| u_{0}
\|^{2}_{H^{s}} + N^{\frac{3}{2}-2s} \lambda^{1-2s} \| u_{0}
\|^{2}_{H^{s}}
\end{array}
\end{equation}
Hence

\begin{equation}
\begin{array}{ll}
\| I u_{\lambda}(0) \|^{4}_{L^{4}} & \lesssim N^{2(1-s)}
\lambda^{1-2s} \| u_{0} \|^{4}_{H^{s}}
\end{array}
\label{Eqn:IuSc}
\end{equation}
By (\ref{Eqn:DIuSc}), (\ref{Eqn:DtIuSc})  and (\ref{Eqn:IuSc}) we
see that there exists $ C_{0}=C_{0} \left( \| u_{0} \|_{H^{s}}, \|
u_{1} \|_{H^{s-1}} \right)$ such that if $\lambda$ satisfies

\begin{equation}
\begin{array}{ll}
\lambda & =  C_{0} N^{\frac{{2(1-s)}}{2s-1}}
\end{array}
\label{Eqn:UpperBdLambda}
\end{equation}
then

\begin{equation}
\begin{array}{ll}
E \left( I u_{\lambda}(0) \right)  & \leq  \frac{1}{2}
\end{array}
\label{Eqn:InitTruncNrjEst}
\end{equation}

\item \textbf{Boundedness of the mollified energy}. Let $F_{T}$ denote the following set

\begin{equation}
\begin{array}{ll}
F_{T}= \left\{ T^{'} \in [0, \,T]: \sup_{t \in [0, \, \lambda
T^{'}]} E \left( I u_{\lambda}(t) \right) \leq 1 \, \mathrm{and} \,
\| I u_{\lambda} \|_{L_{t}^{6} \left( [0, \, \lambda T^{'}] \right)
L_{x}^{6}} \leq (16 C^{2}_{s})^{\frac{1}{6}} +1 \right\}
\end{array}
\nonumber
\end{equation}
with $C_{s}$ being the constant determined by $\lesssim$ in
(\ref{Eqn:SobolevRad}) and $\lambda$ satisfying
(\ref{Eqn:UpperBdLambda}). We claim that $F_{T}$ is the whole set
$[0, \, T]$ for $N=N(T) >> 1$ to be chosen later. Indeed

\begin{itemize}

\item $F_{T} \neq \emptyset$ since $0 \in F_{T}$ by (\ref{Eqn:InitTruncNrjEst}).

\item $F_{T}$ is closed by continuity and by the dominated convergence theorem

\item $F_{T}$ is open. Let $\widetilde{T^{'}} \in F_{T}$. By continuity there exists $\delta > 0$ such that
for every $T^{'} \in \left( \widetilde{T^{'}} - \delta, \,
\widetilde{T^{'}} + \delta \right) \cap [0, \, T]$ we have

\begin{equation}
\begin{array}{ll}
\sup_{t \in [0, \, \lambda T^{'}]} E \left( I u_{\lambda} (t)
\right) & \leq 2
\end{array}
\label{Eqn:HypInducNrj}
\end{equation}
and

\begin{eqnarray}
\| I u_{\lambda} \|_{L_{t}^{6} \left( [0, \, \lambda T^{'}] \right)
L_{x}^{6}} \leq  (16 C_{s}^{2})^{\frac{1}{6}} + 2
\label{Eqn:HypInducLgEst}
\end{eqnarray}
We are interested in generating a partition $\left\{ J_{j} \right\}$
of $[0, \, \lambda T^{'}]$ such that (\ref{Eqn:Condition1}) and
(\ref{Eqn:Condition2}) are satisfied for all $J_{j}$. We describe
now the algorithm.

\emph{Description of the algorithm}. Let $\mathcal{L}$ be the
present list of intervals. Let $L$ be the sum of the lengths of the
intervals making up $\mathcal{L}$. Let $n$ be the number of the last
interval of $\mathcal{L}$. Initially there is no interval and we
start from the time $t=0$. Therefore $\mathcal{L}$ is empty and we
assign the value $0$ to $L$ and $n$. Then as long as $L < \lambda
T^{'}$ do the following

\begin{enumerate}

\item consider $f_{L}(\tau)= \| I u_{\lambda} \|_{L_{t}^{6} \left( [L, \, L + \tau] \right) L_{x}^{6}} - \frac{C_{1}}{\tau^{\frac{1}{3}}}$, $\tau \geq 0 $
with $C_{1}$ defined  in  (\ref{Eqn:Condition1}).

\item since  $f_{L}$ is continuous, does not decrease and $ f_{L}(\tau) \rightarrow  -\infty$ as $\tau \rightarrow 0$, $\tau \geq 0$
there are two options

\begin{itemize}

\item $f_{L}$ is always negative on $[0, \, \lambda T^{'}-L]$: in this case if (\ref{Eqn:Condition2}) is satisfied by
$[L, \, \lambda T^{'}]$ then let $J_{n}:=[L, \, \lambda T^{'}]$. If
not let $J_{n}:=[L, \, L + C_{2} N^{\frac{1-s}{s-\frac{1}{2}}}]$.

\item $f_{L}$ has one and only one root on $[0, \, \lambda T^{'} -L]$: in this case let $\tau_{0}$ be this root.
If (\ref{Eqn:Condition2}) is satisfied by $[L, \, L+ \tau_{0}]$ then
let $J_{n}:= [L, L + \tau_{0}]$. If not let $J_{n} := [L, \, L +
C_{2} N^{\frac{1-s}{s-\frac{1}{2}}}]$.

\end{itemize}

\item assign the value $L + |J_{n}|$ to $L$.

\item assign the value $n+1$ to the variable $n$

\item insert $J_{n}$ into $\mathcal{L}$ so that $\mathcal{L}= \left( J_{j} \right)_{j \in \{1, \, ..., \, n\}}$

\end{enumerate}

\end{itemize}
When we apply this algorithm it is not difficult to see that

\begin{itemize}

\item  $\| I u_{\lambda} \|_{L_{t}^{6} (J_{j}) L_{x}^{6}} = \frac{C_{1}}{|J_{j}|^{\frac{1}{3}}}$ or
$|J_{j}|= C_{2} N^{\frac{1-s}{s - \frac{1}{2}}}$  for every $j \in
\{1, \, ..., \, \card(\mathcal{L}) - 1 \} $

\item $J_{j} \cap J_{k} = \emptyset$ for every $(j,k) \in \{1, \, ..., \, \card(\mathcal{L}) \}^{2}$ such that $j \neq k$

\item $\bigcup_{j=1}^{\card \left( \mathcal{L} \right)} J_{j}$ is a left-closed interval with left endpoint $0$ and included in $[0,\lambda T^{'}]$. Moreover
$\bigcup_{j=1}^{\card \left( \mathcal{L} \right)} J_{j} = [0, \,
\lambda T^{'}]$ if the process is finite.

\end{itemize}
Let

\begin{equation}
\begin{array}{ll}
\mathcal{L}_{1} & = \left\{ J_{j}, \, J_{j} \in \mathcal{L}, \, \| I
u \|_{L_{t}^{6} (J_{j}) L_{x}^{6}} =
\frac{C_{1}}{|J_{j}|^{\frac{1}{3}}} \right\}
\end{array}
\end{equation}
and

\begin{equation}
\begin{array}{ll}
\mathcal{L}_{2} & = \left\{ J_{j}, \, J_{j} \in \mathcal{L}, \,
|J_{j}|=C_{2} N^{\frac{1-s}{s-\frac{1}{2}}} \right\}
\end{array}
\label{Eqn:Boundm1}
\end{equation}
We have $(J_{j})_{j \in \{1, \, ..., \, \card(\mathcal{L})-1 \}}
\subset \mathcal{L}_{1} \cup \mathcal{L}_{2}$. We claim that $
\card(\mathcal{L}_{i}) < \infty $, $i=1,2$. If not let us consider
the $m_{1}$, $m_{2}$ first elements of $\mathcal{L}_{1}$,
$\mathcal{L}_{2}$ respectively. Then

\begin{equation}
\begin{array}{ll}
m_{1} C_{2} N^{\frac{1-s}{s-\frac{1}{2}}}   & \leq \lambda T^{'}
\end{array}
\label{Eqn:Evalm1}
\end{equation}
By H\"older inequality and by (\ref{Eqn:HypInducLgEst}) we have

\begin{equation}
\begin{array}{ll}
m_{2} & = \sum_{j=1}^{m_{2}} |J_{j}|^{-\frac{2}{3}} |J_{j}|^{\frac{2}{3}} \\
& \leq \left( \sum_{j=1}^{m_{2}} \frac{1}{|J_{j}|^{2}}
\right)^{\frac{1}{3}}
\left( \sum_{j=1}^{m_{2}} |J_{j}| \right)^{\frac{2}{3}} \\
& \leq \| Iu \|^{2}_{L_{t}^{6} \left( [0, \, \lambda T^{'}] \right) L_{x}^{6}} \left( \lambda T \right)^{\frac{2}{3}} \\
& \lesssim \left( \lambda T \right)^{\frac{2}{3}}
\end{array}
\label{Eqn:Evalm2}
\end{equation}
Letting $m_{1}$ and $m_{2}$ go to infinity in (\ref{Eqn:Evalm1}) and
(\ref{Eqn:Evalm2}) we have a contradiction. Therefore $\card
(\mathcal{L}) < \infty $ and $\bigcup_{j=1}^{\card \left(
\mathcal{L} \right)} J_{j} = [0, \, \lambda T^{'}]$. Moreover we
have by (\ref{Eqn:UpperBdLambda}), (\ref{Eqn:Evalm1}) and
(\ref{Eqn:Evalm2})

\begin{equation}
\begin{array}{ll}
\card(\mathcal{L}) & \lesssim \left( \lambda T \right)^{\frac{2}{3}} + \frac{\lambda T}{N^{\frac{1-s}{s- \frac{1}{2}}}} + 1 \\
& \lesssim N^{\frac{4(1-s)}{6s-3}} T^{\frac{2}{3}}  + T + 1
\end{array}
\label{Eqn:BoundCardL}
\end{equation}
Now by (\ref{Eqn:InitTruncNrjEst}), (\ref{Eqn:HypInducNrj}),
(\ref{Eqn:BoundCardL}), proposition \ref{prop:LocalBdRad},
\ref{prop:EstNrjRad}, \ref{prop:AlmMorawetz} and
\ref{prop:RemainderMor} we get after iterating

\begin{equation}
\begin{array}{ll}
 \sup_{t \in [0, \, \lambda T^{'}]} E (I u_{\lambda}(t))  - \frac{1}{2}  &
\lesssim \frac{N^{\frac{4(1-s)}{6s-3}} T^{\frac{2}{3}}  + T + 1
}{N^{1-}}
\end{array}
\label{Eqn:EstNrjRadIt}
\end{equation}
and

\begin{equation}
\begin{array}{ll}
\int_{0}^{\lambda T^{'}} \int_{\mathbb{R}^{3}}
\frac{|Iu_{\lambda}(t,x)|^{4}}{|x|} \, dxdt - 2 \left( E (
Iu_{\lambda}(\lambda T^{'}) ) + E(Iu_{\lambda}(0))  \right)  &
\lesssim
\sum_{i=1}^{2} \sum_{j=1}^{\card(\mathcal{L}_{i})} R_{i}(J_{j}) \\
& \lesssim  \frac{N^{\frac{4(1-s)}{6s-3}} T^{\frac{2}{3}}  + T + 1
}{N^{1-}}
\end{array}
\label{Eqn:EstLgEstRadIt}
\end{equation}
By (\ref{Eqn:SobolevRad}), (\ref{Eqn:HypInducNrj}),
(\ref{Eqn:EstLgEstRadIt}) and the inequality $ (1+x)^{\frac{1}{6}}
\leq 1+x$, $x \geq 0$

\begin{equation}
\begin{array}{ll}
\| I u_{\lambda} \|_{L_{t}^{6} ( [0, \, \lambda T^{'}]) L_{x}^{6}} -
(16 C_{s}^{2})^{\frac{1}{6}} & \lesssim \frac{
N^{\frac{4(1-s)}{6s-3}} T^{\frac{2}{3}} + T + 1 }{N^{1-}}
\end{array}
\label{Eqn:EstLgEstRadIt2}
\end{equation}
Let $C^{'}$,$C^{''}$  be the constant determined by $\lesssim$ in
(\ref{Eqn:EstNrjRadIt}), (\ref{Eqn:EstLgEstRadIt2}) respectively.
Since $s>\frac{7}{10}$ we can always choose for every $T>0$ a
$N=N(T) >> 1$ such that

\begin{equation}
\begin{array}{ll}
\frac{ \max{( C^{'}, C^{''} )} N^{\frac{4(1-s)}{6s-3}}
T^{\frac{2}{3}} }{N^{1-}} \leq \frac{1}{6}
\end{array}
\label{Eqn:ChoiceN1}
\end{equation}

\begin{equation}
\begin{array}{ll}
\frac{\max{( C^{'}, C^{''} )} T }{N^{1-}} \leq \frac{1}{6}
\end{array}
\label{Eqn:ChoiceN2}
\end{equation}
and

\begin{equation}
\begin{array}{ll}
\frac{\max{(C^{'}, \, C^{''})}}{N^{1-}} \leq \frac{1}{6}
\end{array}
\label{Eqn:ChoiceN3}
\end{equation}

By (\ref{Eqn:EstNrjRadIt}), (\ref{Eqn:EstLgEstRadIt2}),
(\ref{Eqn:ChoiceN1}), (\ref{Eqn:ChoiceN2}) and (\ref{Eqn:ChoiceN3})
we have $\sup_{t \in [0, \, \lambda T^{'}]} E(I u_{\lambda}(t)) \leq
1 $ and $\| I u_{\lambda} \|_{L_{t}^{6} \left( [0,\, \lambda T^{'}]
\right) L_{x}^{6}} \leq (16 C_{s}^{2})^{\frac{1}{6}} + 1$.

Hence $F_{T}=[0, \, T]$ with $N=N(T)$ satisfying
(\ref{Eqn:ChoiceN1}), (\ref{Eqn:ChoiceN2}) and (\ref{Eqn:ChoiceN3}).

\item \textbf{Conclusion}. Following the $I$- method described
in \cite{almckstt}

\begin{equation}
\begin{array}{ll}
\sup_{t \in [0, \, T]} E \left( I u(t) \right) & = \lambda  \sup_{t \in [0, \, \lambda T]} E \left( (I u)_{\lambda}(t) \right) \\
& \lesssim \lambda \sup_{t \in [0, \, \lambda T]} E (I u_{\lambda}(t)) \\
& \lesssim \lambda
\end{array}
\label{Eqn:ComparNrjScaleRad}
\end{equation}
Combining (\ref{Eqn:ComparNrjScaleRad}) and proposition
\ref{Prop:NrjEst} we have global well-posedness.

Now let $T$ be large. If $ \frac{5}{6} \geq s > \frac{7}{10}$ then
let $N$ such that

\begin{equation}
\begin{array}{ccc}
\frac{0.9}{6} \leq & \frac{ \max{( C^{'}, C^{''} )}
N^{\frac{4(1-s)}{6s-3}} T^{\frac{2}{3}} } {N^{1-}} & \leq
\frac{1}{6}
\end{array}
\label{Eqn:CondNRad1}
\end{equation}
Notice that (\ref{Eqn:ChoiceN2}) and (\ref{Eqn:ChoiceN3}) are also
satisfied. We plug (\ref{Eqn:CondNRad1}) into
(\ref{Eqn:ComparNrjScaleRad}) and we apply proposition
\ref{Prop:NrjEst} to get (\ref{Eqn:LgEstRaduT1}). If $1 > s >
\frac{5}{6}$ then let $N$ such that

\begin{equation}
\begin{array}{ccc}
\frac{0.9}{6} \leq & \frac{ \max{( C^{'}, C^{''} )} T } {N^{1-}} &
\leq \frac{1}{6}
\end{array}
\label{Eqn:CondNRad2}
\end{equation}
Notice that (\ref{Eqn:ChoiceN1}) and (\ref{Eqn:ChoiceN3}) are also
satisfied. We plug (\ref{Eqn:CondNRad2}) into
(\ref{Eqn:ComparNrjScaleRad}) and we apply proposition
\ref{Prop:NrjEst} to get (\ref{Eqn:LgEstRaduT2}).

\end{enumerate}

\section{Proof of the $H^{s}$ norms and mollified energy estimates}
\label{sec:RelHsNrj}

In this section we are interested in proving proposition
\ref{Prop:NrjEst}. By Plancherel theorem

\begin{equation}
\begin{array}{ll}
\| u(T) \|^{2}_{H^{s}} & \lesssim  \| P_{\leq 1} u(T) \|^{2}_{H^{s}}
+ \int_{1\leq |\xi| \leq 2N} |\xi|^{2s} |\hat{u}(T,\xi)|^{2} d\xi +
\int_{|\xi| \geq 2N}|\xi|^{2s} \left| \hat{u}(T,\xi) \right|^{2} \, d\xi \nonumber \\
\end{array}
\end{equation}
But

\begin{equation}
\begin{array}{ll}
\int_{1 \leq |\xi| \leq 2N} |\xi|^{2s} \left| \hat{u}(T,\xi)
\right|^{2} \, d\xi & \leq \int_{|\xi| \leq 2N} |\xi|^{2} \left|
\hat{u}(T,\xi) \right|^{2} \, d\xi \\
& \lesssim  \int_{\mathbb{R}^{3}} \left| D I u (T,x) \right|^{2} \, dx \\
& \lesssim  E \left( Iu(T) \right)
\end{array}
\label{Eqn:HighFrequEIu1}
\end{equation}

\begin{equation}
\begin{array}{ll}
\int_{|\xi| \geq 2N} |\xi|^{2s} |\hat{u}(T,\xi)|^{2} \, d\xi & \leq
\int_{|\xi| \geq 2N} |\xi|^{2} \frac{N^{2(1-s)}}{|\xi|^{2(1-s)}}
\left| \hat{u}(T,\xi)
\right|^{2} \, d\xi  \\
& \lesssim  \int_{\mathbb{R}^{3}} \left| D I u (T,x) \right|^{2} \, dx \\
& \lesssim  E \left( Iu(T) \right)
\end{array}
\label{Eqn:HighFreqEIu2}
\end{equation}
and by the fundamental theorem of calculus and Minkowski inequality

\begin{equation}
\begin{array}{ll}
\| P_{\leq 1} u(T) \|_{H^{s}} & \lesssim \| P_{\leq 1}  u_{0}
\|_{H^{s}} + \int_{0}^{T} \| P_{\leq 1}
\partial_{t} u(t) \|_{H^{s}} \, dt \\
& \lesssim \| u_{0} \|_{H^{s}} + T  \sup_{t \in [0, \, T]} \|
\partial_{t} I u(t) \|_{L^{2}}
\end{array}
\end{equation}
which implies that

\begin{equation}
\begin{array}{ll}
\| P_{\leq 1} u(T) \|^{2}_{H^{s}} & \lesssim \| u_{0} \|^{2}_{H^{s}}
+ T^{2} \sup_{t \in [0, \, T]} E \left( Iu(t) \right)
\end{array}
\label{Eqn:LowFrequEIu}
\end{equation}
We also have

\begin{equation}
\begin{array}{ll}
\| \partial_{t} u(T) \|^{2}_{H^{s-1}} & \lesssim  E \left( Iu(T)
\right)
\end{array}
\label{Eqn:DerivuEIu}
\end{equation}
Combining (\ref{Eqn:HighFrequEIu1}),
(\ref{Eqn:HighFreqEIu2}),(\ref{Eqn:LowFrequEIu}) and
(\ref{Eqn:DerivuEIu}) we get (\ref{Eqn:NrjEstim1}).

\section{Proof of the local boundedness estimate}
\label{sec:LocalBd}

We are interested in proving proposition \ref{prop:LocalBdRad}  in
this section. In what follows we also assume that $J=[0, \, \tau]$:
the reader can check after reading the proof that the other cases
can be reduced to that one.

Before starting the proof let us state the following lemma

\begin{lem}{\textbf{"Strichartz estimates with derivative"}}
Let $m \in [0, \, 1]$ and $ 0 \leq \tau < \infty$. If $u$ satisfies
the IVP problem

\begin{equation}
\left\{
\begin{array}{ccc}
\Box u & = & F  \\
u(t=0) & = & f   \\
\partial_{t} u(t=0) & = & g
\end{array}
\right.
\label{Eqn:EqWaveGen}
\end{equation}
then we have the $m$- Strichartz estimate with derivative

\begin{equation}
\begin{array}{ll}
\| u \|_{L_{t}^{q}\left( [0, \, \tau] \right)  L_{x}^{r}} + \|
\partial_{t} D^{-1} u \|_{L_{t}^{q} \left( [0, \, \tau] \right)
L_{x}^{r}} & \lesssim  \| f \|_{\dot{H}^{m}} + \| g
\|_{\dot{H}^{m-1}} + \| F \|_{L_{t}^{\tilde{q}} \left( [0, \, \tau]
\right) L_{x}^{\tilde{r}}}
\end{array}
\label{Eqn:StrEstDeriv}
\end{equation}
for $(q,r) \in \mathcal{W} $, $(\tilde{q},\tilde{r}) \in
\widetilde{\mathcal{W}}$ and $(q,r,\tilde{q},\tilde{r})$ satisfying
the gap condition

\begin{equation}
\begin{array} {lll}
\frac{1}{q} + \frac{3}{r} & =  \frac{3}{2} -m  & =
\frac{1}{\tilde{q}} + \frac{3}{\tilde{r}} -2
\end{array}
\end{equation}
\label{lem:StrEstDeriv}
\end{lem}
We postpone the proof of lemma \ref{lem:StrEstDeriv} to subsection
\ref{subsec:PfLemStrDer}. Assuming that is true we now show how
lemma \ref{lem:StrEstDeriv} implies proposition
\ref{prop:LocalBdRad}.

Multiplying the $m$-Strichartz estimate with derivative
(\ref{Eqn:StrEstDeriv}) by $D^{1-m}I$ we get

\begin{eqnarray}
 Z_{m,s}(\tau) & \lesssim & \| D I u_{0} \|_{L^{2}} + \| I u_{1} \|_{L^{2}} + \| D^{1-m} I F \|
 _{L_{t}^{\tilde{q}}([0..\tau]) L_{x}^{\tilde{r}} } \nonumber \\
 & \lesssim & 1 + \| D^{1-m} I F \|_{L_{t}^{\tilde{q}}([0, \,\tau]) L_{x}^{\tilde{r}} }
 \label{Eqn:ZEst}
\end{eqnarray}
The remainder of proof is divided into three steps.

\begin{itemize}

\item \textbf{First Step}
First we assume that $m \leq s$. Notice that the point
$(\frac{1}{1-s},6)$ is $s$-wave admissible. In this case we get from
the fractional Leibnitz rule the H\"older in time and the H\"older
in space inequalities

\begin{equation}
\begin{array}{ll}
Z_{m,s}(\tau) & \lesssim  1 + \| D^{1-m} I (uuu) \|_{L_{t}^{1}([0, \, \tau]) L_{x}^{\frac{6}{5-2m}}} \\
& \lesssim  1+ \| D^{1-m} I u \|_{L_{t}^{\infty}([0, \, \tau])
L_{x}^{\frac{6}{3-2m}}} \| u \|_{L_{t}^{2} \left( [0, \, \tau]
\right)
L_{x}^{6}}^{2} \\
& \lesssim  1+ Z_{m,s} (\tau) \left( \tau^{\frac{1}{3}} \| P_{\leq
N} u \|_{L_{t}^{6} \Ivl L_{x}^{6}} + \tau^{s-\frac{1}{2}} \| P_{>N}
u \|_{L_{t}^{\frac{1}{1-s}} \Ivl
L_{x}^{6}}  \right)^{2} \\
& \lesssim  1 + Z_{m,s}(\tau) \left( \tau^{\frac{1}{3}} \| Iu
\|_{L_{t}^{6} \Ivl  L_{x}^{6}} + \tau^{s-\frac{1}{2}} \frac{ \|
D^{1-s} Iu  \|_{L_{t}^{\frac{1}{1-s}}\left( [0, \, \tau] \right)
L_{x}^{6}}}{N^{1-s}}  \right)^{2}  \\
& \lesssim  1 + Z_{m,s}(\tau) \left( \tau^{\frac{1}{3}} \| Iu
\|_{L_{t}^{6} \Ivl L_{x}^{6}} + \tau^{s-\frac{1}{2}}
\frac{Z_{s,s}(\tau)}{N^{1-s}}  \right)^{2} \label{Eqn:IndZms}
\end{array}
\end{equation}
Assume $m=s$. Then if we apply a continuity argument to
(\ref{Eqn:IndZms}) we get from the inequalities
(\ref{Eqn:Condition1}) and (\ref{Eqn:Condition2})

\begin{equation}
\begin{array}{ll}
Z_{s,s}(\tau) & \lesssim  1
\end{array}
\label{Eqn:BdZss}
\end{equation}
Now assume  $m<s$. Then if we apply a continuity argument to
(\ref{Eqn:IndZms}) and the inequalities (\ref{Eqn:Condition1}) and
(\ref{Eqn:BdZss})  we get

\begin{equation}
\begin{array}{ll}
Z_{m,s}(\tau) & \lesssim  1
\end{array}
\label{Eqn:BdZms}
\end{equation}

\item \textbf{Second Step}
We assume $m>s$. By (\ref{Eqn:IndZms}), (\ref{Eqn:BdZss}),
(\ref{Eqn:BdZms}), (\ref{Eqn:Condition1}) and (\ref{Eqn:Condition2})
we have

\begin{equation}
\begin{array}{ll}
\| D^{1-r} I (uuu) \|_{L_{t}^{1} \Ivl L_{x}^{\frac{6}{5-2r}}} &
\lesssim  Z_{r,s}(\tau) \left( \tau^{\frac{1}{3}} \| Iu
\|_{L_{t}^{6} \Ivl L_{x}^{6}} +
\frac{ \tau^{s-\frac{1}{2}} Z_{s,s}(\tau)}{N^{1-s}}  \right)^{2} \\
& \lesssim  1
\end{array}
\label{Eqn:BdFcTerm}
\end{equation}
for $r \leq s$. The inequality

\begin{equation}
\begin{array}{ll}
\| D^{1-m} I (uuu) \|_{L_{t}^{1} \Ivl L_{x}^{\frac{6}{5-2m}}} &
\lesssim  \| D^{1-r} I (uuu) \|_{L_{t}^{1} \Ivl
L_{x}^{\frac{6}{5-2r}}}
\end{array}
\label{Eqn:SobIneqm}
\end{equation}
follows from the application of Sobolev  homogeneous embedding. We
get from (\ref{Eqn:ZEst}), (\ref{Eqn:BdFcTerm}) and
(\ref{Eqn:SobIneqm})

\begin{equation}
\begin{array}{ll}
Z_{m,s}(\tau) & \lesssim  1 + \| D^{1-m} I (uuu) \|_{L_{t}^{1} \left( [ 0, \, \tau] \right)  L_{x}^{\frac{6}{5-2m}}} \\
& \lesssim  1 + \| D^{1-r} I (uuu) \|_{L_{t}^{1} \left( [ 0, \, \tau] \right)  L_{x}^{\frac{6}{5-2r}}}  \\
& \lesssim  1
\end{array}
\end{equation}

\end{itemize}

\subsection{Proof of lemma \ref{lem:StrEstDeriv}}
\label{subsec:PfLemStrDer}

By decomposition it suffices to prove that $u_{l}^{1}(t)=e^{ \pm
itD}f$, $u_{l}^{2}(t)=\frac{e^{\pm itD}}{D}g$ and $u_{n}(t)=
\int_{0}^{t} D^{-1} \sin { \left( (t-t^{'})D \right) }  F \, dt^{'}$
satisfy (\ref{Eqn:StrEstDeriv}).

We have $ \partial_{t} u_{l}^{1}(t)= \pm i D e^{ \pm itD} f $ and
$\partial_{t} u_{l}^{2} = \pm e^{\pm it D} g$. We know from the
Strichartz estimates that

\begin{equation}
\begin{array}{ll}
\| D^{-1} \partial_{t} u_{l}^{1} \|_{L_{t}^{q} \Ivl L_{x}^{r} } &
\lesssim
\| e^{\pm it D} f \|_{L_{t}^{q} \Ivl  L_{x}^{r}} \\
& \lesssim  \| f \|_{\dot{H}^{m}}
\end{array}
\label{Eqn:Lin1}
\end{equation}
and

\begin{equation}
\begin{array}{ll}
\| D^{-1} \partial_{t} u_{l}^{2} \|_{L_{t}^{q} \Ivl L_{x}^{r}} & =
\| e^{\pm itD} D^{-1} g \|_{L_{t}^{q} \Ivl L_{x}^{r}}  \\
& \lesssim  \| D^{-1} g \|_{\dot{H}^{m}} \\
& \lesssim  \| g \|_{\dot{H}^{m-1}}
\end{array}
\label{Eqn:Lin2}
\end{equation}
We also have

\begin{equation}
\begin{array}{ll}
D^{-1} \partial_{t} u_{n}(t) & =   \int_{0}^{t} \cos{ \left(
(t-t^{'})D \right) } F(t{'}) \, dt^{'}
\end{array}
\end{equation}
and by the Strichartz  estimates

\begin{equation}
\begin{array}{ll}
\| D^{-1} \partial_{t} u_{n} \|_{L_{t}^{q} \Ivl L_{x}^{r}} &
\lesssim
\|  \int_{0}^{t} D^{-1} e^{ i(t-t^{'}) D } F(t^{'}) dt^{'}  \|_{L_{t}^{q} \Ivl L_{x}^{r}}  \\
&  + \|  \int_{0}^{t} D^{-1} e^{ -i(t-t^{'}) D } F(t^{'}) dt^{'}  \|_{L_{t}^{q} \Ivl L_{x}^{r}} \\
& \lesssim  \| F \|_{L_{t}^{\tilde{q}} \Ivl L_{x}^{\tilde{r}}}
\end{array}
\label{Eqn:ResNlin}
\end{equation}
(\ref{Eqn:StrEstDeriv}) follows from (\ref{Eqn:Lin1}),
(\ref{Eqn:Lin2}) and (\ref{Eqn:ResNlin}).

\section{Proof of almost conservation law}
\label{sec:AlmCon}

Now we prove proposition \ref{prop:EstNrjRad}. In what follows we
also assume that $J=[0, \, \tau]$: the reader can check after
reading the proof that the other cases can be reduced to that one.

Let $\tau_{0} \in J$. It suffices to prove

\begin{equation}
\begin{array}{ll}
\left| E \left( Iu(\tau_{0}) \right) - E \left( Iu(0) \right)
\right| & \lesssim \frac{Z^{4}(\tau)}{N^{1-}}
\end{array}
\end{equation}

In what follows we also assume that $\tau_{0}= \tau$: the reader can
check after reading the proof that the other cases can be reduced to
this one.

The Plancherel formula and the fundamental theorem of calculus yield

\begin{equation}
\begin{array}{ll}
 E \left( Iu(\tau) \right) - E \left( Iu(0) \right) & =
 \int_{0}^{\tau} \int_{\xi_{1} + ... + \xi_{4} = 0}
\mu(\xi_{2},\xi_{3},\xi_{4}) \widehat{\partial_{t} I u}(t,\xi_{1})
\widehat{I u}(t,\xi_{2}) \widehat{I u} (t,\xi_{3}) \widehat{I u}
(t,\xi_{4}) \, d\xi_{2} d\xi_{3} d \xi_{4} dt
\end{array}
\nonumber
\end{equation}
with

\begin{equation}
\begin{array}{ll}
\mu(\xi_{2},\xi_{3},\xi_{4}) & = 1 - \frac{m(\xi_{2}+ \xi_{3} +
\xi_{4})}{m(\xi_{2}) m(\xi_{3}) m(\xi_{4})}
\end{array}
\label{Eqn:DfnMu}
\end{equation}
We are left to prove

\begin{equation}
\begin{array}{ll}
\left| \int_{0}^{\tau} \int_{\xi_{1} + ... + \xi_{4} = 0}
\mu(\xi_{2},\xi_{3},\xi_{4}) \widehat{\partial_{t} I u}(t,\xi_{1})
\hat{I u}(t, \xi_{2}) \widehat{I u }(t,\xi_{3}) \widehat{I u} (t,
\xi_{4}) \, d\xi_{2}d\xi_{3}d \xi_{4} \, dt \right| & \lesssim
\frac{Z^{4}(\tau)}{N^{1-}}
\end{array}
\label{Eqn:AlmConsToProve}
\end{equation}
We perform a Paley-Littlewood decomposition to prove
(\ref{Eqn:AlmConsToProve}). Let $u_{i}=P_{N_{i}} u$ with $i \in \{1,
\, ..., \, 4 \}$ and let

\begin{equation}
\begin{array}{ll}
X & =  \left| \int_{0}^{\tau} \int_{\xi_{1}+...+\xi_{4}=0}
\mu(\xi_{2},\xi_{3},\xi_{4}) \widehat{\partial_{t} I u_{1}}(t,
\xi_{1}) \widehat{I u_{2}}(t,\xi_{2}) \widehat{I u_{3}}(t,\xi_{3})
\widehat{I u_{4}} (t,\xi_{4}) \,  d\xi_{2}d\xi_{3}d \xi_{4}dt
\right|
\end{array}
\label{Eqn:DfnXAlmRad}
\end{equation}
There are different cases resulting from this Paley-Littlewood
analysis and we describe now the strategy to estimate
(\ref{Eqn:AlmConsToProve}). We suggest that the reader at first
ignores the second and third steps of the description and the
$N_{j}^{\pm}$ appearing in the study of these cases to solve the
summation issue. \vspace{5 mm}

\emph{Description of the strategy}
\begin{enumerate}

\item We follow \cite{morckstt} to estimate $X$. First we recall the following Coifman-Meyer theorem
\cite{coifmey2}, p179 for a class of multilinear operators

\begin{thm}  {\textbf{''Coifman Meyer multiplier theorem''} }
Consider an infinitely differentiable symbol $\sigma:
\mathbb{R}^{nk} \rightarrow \mathbb{C}$ so that for all $\alpha \in
\mathbb{N}^{nk}$ there exists $c(\alpha)$ such that for all
 $ \xi= \left( \xi_{1},...,\xi_{k} \right) \in \mathbb{R}^{nk}$

\begin{eqnarray}
\left| \partial_{\xi}^{\alpha} \sigma (\xi) \right| & \leq
\frac{c(\alpha)} { \left( 1+ |\xi| \right)^{| \alpha|}}
\label{Eqn:MeyMultBound}
\end{eqnarray}
Let $\Lambda_{\sigma}$ be the multilinear operator

\begin{equation}
\begin{array}{ll}
\Lambda_{\sigma} (f_{1},...,f_{k})(x) & =  \int_{\mathbb{R}^{nk}}
e^{ix \cdot (\xi_{1}+...+\xi_{k})} \sigma(\xi_{1},...,\xi_{k})
\widehat{f_{1}}(\xi_{1})...\widehat{f_{k}}(\xi_{k}) \, d \xi_{1}...d
\xi_{k}
\end{array}
\label{Eqn:ThmCoifMey}
\end{equation}
Assume that $q_{j} \in (1,\infty)$, $j \in \{1, \, ..., \, k \}$ are
such that $\frac{1}{q}=\frac{1}{q_{1}}+...+\frac{1}{q_{k}} \leq 1$.
Then there is a constant $C=C \left( q_{j},n,k,c(\alpha) \right)$ so
that for all Schwarz class functions $f_{1},...,f_{k}$

\begin{eqnarray}
\| \Lambda_{\sigma}(f_{1},...,f_{k}) \|_{L^{q}(\mathbb{R}^{n})} &
\leq & C \| f_{1} \|_{L^{q_{1}} (\mathbb{R}^{n})}...\| f_{k}
\|_{L^{q_{k}} (\mathbb{R}^{n})}
\end{eqnarray}
\label{thm:CoifMeyer}
\end{thm}
Then we proceed as follows. We seek a pointwise bound on the symbol

\begin{equation}
\begin{array}{ll}
\left| \mu(\xi_{2},\xi_{3},\xi_{4})  \right| & \leq  B \left(
N_{2},N_{3},N_{4} \right)
\end{array}
\label{Eqn:DfnB}
\end{equation}
We factor $B=B(N_{2},N_{3},N_{4})$ out of the right side of
(\ref{Eqn:DfnXAlmRad}) and we are left to evaluate

\begin{equation}
\begin{array}{l}
B  \int_{0}^{\tau} \int_{\mathbb{R}^{3}}
\widehat{\Lambda_{\frac{\mu}{B}} (
\partial_{t} I u_{1}(t), I u_{2}(t), I u_{3}(t) )}
(\xi_{4}) \widehat {I u_{4}}(t,\xi_{4}) \, d \xi_{4} \, dt
\end{array}
\nonumber
\end{equation}
We notice that the multiplier $\frac{\mu}{B}$ satisfy the bound
(\ref{Eqn:MeyMultBound}) and by the Plancherel  theorem, H\"older
inequality, theorem \ref{thm:CoifMeyer} and Bernstein inequalities
we have

\begin{equation}
\begin{array}{ll}
X & \lesssim  B \| \partial_{t} I u_{1} \|_{L_{t}^{p_{1}} ([0, \,
\tau]) L_{x}^{q_{1}} }
\| I u_{2} \|_{L_{t}^{p_{2}} ([0, \, \tau]) L_{x}^{q_{2}}}... \| I u_{4} \|_{L_{t}^{p_{4}} ([0, \, \tau]) L_{x}^{q_{4}}}  \\
& \lesssim  B  N_{1}^{m_{1}} N_{2}^{m_{2}-1}...N_{4}^{m_{4}-1} \|
\partial_{t} D^{-m_{1}} I u_{1} \|_{L_{t}^{p_{1}} ([0, \, \tau])
L_{x}^{q_{1}} } \| D^{1-m_{2}} I u_{2} \|_{L_{t}^{p_{2}} ([0, \, \tau]) L_{x}^{q_{2}}}... \| D^{1-m_{4}} I u_{4} \|_{L_{t}^{p_{4}} ([0, \, \tau]) L_{x}^{q_{4}}} \\
&  \lesssim B  N_{1}^{m_{1}} N_{2}^{m_{2}-1}...N_{4}^{m_{4}-1}
Z^{4}(\tau)
\end{array}
\label{Eqn:AlmconsRadpj}
\end{equation}
with $(p_{j},q_{j})$ such that $p_{j} \in [1, \, \infty]$ and $q_{j}
\in (1, \, \infty)$ for $j=\{1, \, ..., \, 4 \}$, $\sum_{j=1}^{4}
\frac{1}{p_{j}} = 1$, $\sum_{j=1}^{4} \frac{1}{q_{j}} =1$,
$(p_{j},q_{j})$ $m_{j}$-wave admissible for some $m_{j}^{'} \, s$
such that $0 \leq m_{j} < 1$ and $\frac{1}{p_{j}} + \frac{1}{q_{j}}
=\frac{1}{2}$ \footnote{in other words $(p_{j},q_{j}) = \left(
\frac{2}{m_{j}}, \, \frac{2}{1-m_{j}} \right)$}.

\item The series must be summable. Therefore in some cases we might create $N_{k}^{\pm}$ for some $k^{'}s$
by considering slight variations $(p_{k} \pm , \, q_{k} \pm) \in [1,
\, \infty] \times (1, \, \infty)$ of $(p_{k},q_{k})$ that are $m_{k}
\, \pm$ - wave admissible and such that $\frac{1}{p_{k} \pm} +
\frac{1}{q_{k} \pm} =\frac{1}{2}$. For instance if we create slight
variations $(p_{2}+,q_{2}-)$, $(p_{4}-,q_{4}+)$ of $(p_{2},q_{2})$,
$(p_{4},q_{4})$ respectively we have

\begin{equation}
\begin{array}{ll}
\| I u_{2} \|_{L_{t}^{p_{2}+} L_{x}^{q_{2}-}} & \lesssim N_{2}^{-} N_{2}^{m_{2}-1} \| D^{1-(m_{2}-)} I u_{2} \|_{L_{t}^{p_{2}+} L_{x}^{q_{2}-}} \\
\| I u_{4} \|_{L_{t}^{p_{4}-} L_{x}^{q_{4}+}} & \lesssim N_{4}^{+}
N_{4}^{m_{4}-1} \| D^{1-(m_{4}+)} I u_{4} \|_{L_{t}^{p_{4}-}
L_{x}^{q_{4}+}}
\end{array}
\label{Eqn:Ex1DirectCreat}
\end{equation}
and (\ref{Eqn:AlmconsRadpj}) becomes

\begin{equation}
\begin{array}{ll}
X & \lesssim B N_{2}^{-} N_{4}^{+}  N_{1}^{m_{1}}
N_{2}^{m_{2}-1}...N_{4}^{m_{4}-1} Z^{4}(\tau)
\end{array}
\end{equation}

\item When we deal with low frequencies, i.e $N_{k}<1$ for some $k \in \{1,...,4 \}$ we might consider generating $N_{k}^{+}$ by
creating a variation $(2+, \infty-)$ of $(2,\infty)$. Such a task
cannot be directly performed since we unfortunately have

\begin{equation}
\begin{array}{ll}
\| I u_{k} \|_{L_{t}^{2+} L_{x}^{\infty -}} & \lesssim N_{k}^{-} \| D^{1-(1-)} I u_{k} \|_{L_{t}^{2+} L_{x}^{\infty -}} \\
& \lesssim N_{k}^{-} Z(\tau)
\end{array}
\label{Eqn:Ex2DirectCreat}
\end{equation}
But we can indirectly create $N_{k}^{+}$ by appropriately using
H\"older in time inequality. Indeed if $\epsilon > 0$, $\epsilon^{'}
> 0 $ and $ \epsilon^{''} > 0 $ are such that $\frac{\epsilon}{2}=
\frac{\epsilon^{'}}{2} - \frac{\epsilon^{''}}{3}$ we get from
Bernstein inequalities, H\"older in time inequality and Sobolev
homogeneous embedding

\begin{equation}
\begin{array}{ll}
\| I u_{k} \|_{L_{t}^{\frac{2}{1- \epsilon}}(\Ivl)
L_{x}^{\frac{2}{\epsilon}}} &  \lesssim  N_{k}^{\epsilon^{'}}
 \| D^{-\epsilon^{'}} I u_{k}
\|_{L_{t}^{\frac{2}{1- \epsilon}}(\Ivl)
L_{x}^{\frac{2}{\epsilon}}} \\
& \lesssim N_{k}^{\epsilon^{'}} \tau^{\frac{\epsilon^{'}
-\epsilon}{2}} \| D^{-\epsilon^{'}} I u_{k} \|_{L_{t}^{\frac{2}{1-
\epsilon^{'}}}(\Ivl)
L_{x}^{\frac{2}{\epsilon}}} \\
&  \lesssim N_{k}^{\epsilon^{'}} \tau^{\frac{\epsilon^{'}
-\epsilon}{2}} \| D^{-\epsilon^{'}+ \epsilon^{''}} I u_{k}
\|_{L_{t}^{\frac{2}{1- \epsilon^{'}}}(\Ivl)
L_{x}^{\frac{2}{\epsilon^{'}}}} \\
& \lesssim  N_{k}^{\epsilon^{''} - \epsilon^{'}}
\tau^{\frac{\epsilon^{'} -\epsilon}{2}} \| D^{ 1- (1- \epsilon^{'})}
I u_{k} \|_{L_{t}^{\frac{2}{1- \epsilon^{'}}}(\Ivl)
L_{x}^{\frac{2}{\epsilon^{'}}}} \\
& \lesssim N_{k}^{\epsilon^{''} - \epsilon^{'}}
\tau^{\frac{\epsilon^{'} -\epsilon}{2}} Z(\tau)
\end{array}
\label{Eqn:IndirectStepsQuad}
\end{equation}
We would like $\epsilon^{''} > \epsilon^{'}$. A quick computation
show that it suffices that $ \epsilon^{'} > 3 \epsilon $. Letting
$\epsilon^{'}= 5 \epsilon$ we get

\begin{equation}
\begin{array}{ll}
\| I u_{k} \|_{L_{t}^{\frac{2}{1- \epsilon}}(\Ivl)
L_{x}^{\frac{2}{\epsilon}}} & \lesssim N_{k}^{\epsilon} \tau^{ 2
\epsilon } Z(\tau)
\end{array}
\end{equation}
Now if we choose $\epsilon > 0 $ so small that $ |\tau|^{2 \epsilon}
\leq 2 $ we eventually get

\begin{equation}
\begin{array}{ll}
\| I u_{k} \|_{L_{t}^{2+}(\Ivl) L_{x}^{\infty -}} & \lesssim
N_{k}^{+}  Z(\tau)
\end{array}
\label{Eqn:IndirectRes}
\end{equation}
For the remainder of the paper we say that we directly create
$N_{k}^{\pm}$ if we directly use Bernstein inequality like in
(\ref{Eqn:Ex1DirectCreat}) or (\ref{Eqn:Ex2DirectCreat}) and we say
that we indirectly create $N_{k}^{+}$  if we also use H\"older in
time inequality to get (\ref{Eqn:IndirectRes}).

\end{enumerate}
This ends the general description of the strategy. \vspace{5mm}

Let us get back to the proof. By symmetry we may assume that $N_{2}
\geq N_{3} \geq N_{4}$. There are several cases

\begin{itemize} 

\item \textbf{Case} $\mathbf{1}$: $N >> N_{2} \geq N_{3}$. In this case $X=0$
since $\mu=0$.

\item \textbf{Case} $\mathbf{2}$: $N_{2} \gtrsim N >> N_{3}$

In this case we have
\begin{equation}
\begin{array}{ll}
 \left| \mu(\xi_{2},..,\xi_{4}) \right| & \lesssim \frac{|\nabla m(\xi_{2})| |\xi_{3} + \xi_{4}|}{m(\xi_{2})} \\
 & \lesssim \frac{N_{3}}{N_{2}}
\end{array}
\label{Eqn:EstMultRad1}
\end{equation}
We also get $N_{1} \sim N_{2}$ from the convolution constraint
$\xi_{1}+...+ \xi_{4}=0$.

We assume that $N_{4} \geq 1$. By (\ref{Eqn:EstMultRad1}) and by the
Bernstein inequalities we have

\begin{equation}
\begin{array}{ll}
X & \lesssim \frac{N_{3}}{N_{2}} \| \partial_{t} I u_{1}
\|_{L_{t}^{6-} \Ivl L_{x}^{3+}} \| I u_{2} \|_{L_{t}^{6} \Ivl
L_{x}^{3}}
\| I u_{3} \|_{L_{t}^{6} \Ivl L_{x}^{3}} \| I u_{4} \|_{L_{t}^{2+} \Ivl L_{x}^{\infty -}} \\
& \lesssim N_{1}^{+} N_{4}^{-} \frac{N_{3}}{N_{2}}
N_{1}^{\frac{1}{3}} N_{2}^{-\frac{2}{3}} N_{3}^{-\frac{2}{3}}
\| \partial_{t} D^{- \left( \frac{1}{3}+ \right) } I u_{1} \|_{L_{t}^{6-} \Ivl L_{x}^{3+}} \\
& \| D^{1- \frac{1}{3} } I u_{2} \|_{L_{t}^{6} \Ivl L_{x}^{3}}
\| D^{1-\frac{1}{3}}  I u_{3} \|_{L_{t}^{6} \Ivl L_{x}^{3}}  \| D^{1-(1-)} I u_{4} \|_{L_{t}^{2+} \Ivl L_{x}^{\infty-}} \\
& \lesssim \frac{N_{2}^{-} N_{4}^{-}}{N^{1-}} Z^{4}(\tau)
\end{array}
\nonumber
\end{equation}
after directly creating $N_{1}^{+}$ and $N_{4}^{-}$. If $N_{4} < 1$
the proof is similar except that we indirectly create $N_{4}^{+}$ to
get $X \lesssim \frac{N_{2}^{-} N_{4}^{+}}{N^{1-}} Z^{4}(\tau) $.
This makes the summation possible. We get (\ref{Eqn:AlmConsToProve})
after summation.

\item \textbf{Case} $\mathbf{3}$: $ N_{3} \gtrsim N >> N_{4}$

In this case we have

\begin{equation}
\begin{array}{ll}
\left| \mu(\xi_{2},..,\xi_{4}) \right| & \lesssim
\frac{m(\xi_{1})}{m(\xi_{2}) m(\xi_{3}) m(\xi_{4})}
\end{array}
\label{Eqn:EstMultRad2}
\end{equation}
There are two subcases

\begin{itemize}

\item \textbf{Case} $\mathbf{3.a}$: $N_{1} \sim N_{2}$

We assume that $N_{4} \geq 1$. By (\ref{Eqn:EstMultRad2}) we have

\begin{equation}
\begin{array}{ll}
X & \lesssim \frac{N_{3}^{1-s}}{N^{1-s}} \| \partial_{t} I u_{1}
\|_{L_{t}^{6-} \Ivl L_{x}^{3+}} \| I u_{2} \|_{L_{t}^{6} \Ivl
L_{x}^{3}}
\| I u_{3} \|_{L_{t}^{6} \Ivl L_{x}^{3}} \| I u_{4} \|_{L_{t}^{2+} \Ivl L_{x}^{\infty -}} \\
& \lesssim N_{1}^{+} N_{4}^{-} \frac{N_{3}^{1-s}}{N^{1-s}}
N_{1}^{\frac{1}{3}} N_{2}^{-\frac{2}{3}} N_{3}^{-\frac{2}{3}}   \|
\partial_{t} D^{- \left( \frac{1}{3}+ \right) } I u_{1}
\|_{L_{t}^{6-} \Ivl L_{x}^{3+}}
\| D^{1- \frac{1}{3} } I u_{2} \|_{L_{t}^{6} \Ivl L_{x}^{3}} \\
& \| D^{1-\frac{1}{3}}  I u_{3} \|_{L_{t}^{6} \Ivl L_{x}^{3}}  \| D^{1-(1-)} I u_{4} \|_{L_{t}^{2+} \Ivl L_{x}^{\infty-}} \\
& \lesssim \frac{N_{2}^{-} N_{4}^{-}}{N^{1-}} Z^{4}(\tau)
\end{array}
\nonumber
\end{equation}
after directly creating $N_{1}^{+}$ and $N_{4}^{-}$. If $N_{4} < 1$
the proof is similar except that we indirectly create $N_{4}^{+}$.
We get (\ref{Eqn:AlmConsToProve}) after summation.

\item \textbf{Case} $\mathbf{3.b}$: $N_{1} << N_{2}$

In this case by the convolution constraint $\xi_{1}+...+\xi_{4}=0$
we have $N_{2} \sim N_{3}$. There are two subcases

\begin{itemize}

\item \textbf{Case} $\mathbf{3.b.1}$: $N_{1} << N$

We assume that $N_{1} \geq 1$ and $N_{4} \geq 1$. By
(\ref{Eqn:EstMultRad2}) we have

\begin{equation}
\begin{array}{ll}
X & \lesssim \frac{N_{2}^{2(1-s)}}{N^{2(1-s)}} \| \partial_{t} I
u_{1} \|_{L_{t}^{6-} \Ivl L_{x}^{3+}} \| I u_{2} \|_{L_{t}^{6} \Ivl
L_{x}^{3}}
 \| I u_{3} \|_{L_{t}^{6} \Ivl L_{x}^{3}} \\ &  \| I u_{4} \|_{L_{t}^{2+} \Ivl L_{x}^{\infty -}} \\
& \lesssim  N_{1}^{+} N_{4}^{-} \frac{N_{2}^{2(1-s)}}{N^{2(1-s)}}
N_{1}^{\frac{1}{3}} N_{2}^{-\frac{2}{3}} N_{3}^{-\frac{2}{3}}
\| \partial_{t} D^{- \left( \frac{1}{3}+ \right) } I u_{1} \|_{L_{t}^{6-} \Ivl L_{x}^{3+}} \\
& \| D^{1- \frac{1}{3} } I u_{2} \|_{L_{t}^{6} \Ivl L_{x}^{3}}  \|
D^{1-\frac{1}{3}}  I u_{3} \|_{L_{t}^{6} L_{x}^{3}}
\| D^{1-(1-)} I u_{4} \|_{L_{t}^{2+} \Ivl L_{x}^{\infty-}} \\
& \lesssim \frac{N_{1}^{-} N_{2}^{-}  N_{4}^{-}}{N^{1-}} Z^{4}(\tau)
\end{array}
\nonumber
\end{equation}
after directly creating $N_{1}^{+}$ and $N_{4}^{-}$. If $N_{1} < 1$
and $N_{4} < 1$ the proof is similar except that we indirectly
create $N_{4}^{+}$ and we substitute $N_{1}^{-}$ for $N_{1}^{+}$.
The proof for the other cases \footnote{i.e $N_{1} \geq 1$, $N_{4}
\leq 1$ or $N_{1} \leq 1$, $N_{4} \geq 1$} is a slight variant to
that for the case $N_{1} \geq 1$, $N_{4} \geq 1$ and that for the
case $N_{1} < 1$, $N_{4}<1$. Details are left to the reader. We get
(\ref{Eqn:AlmConsToProve}) after summation.

\item \textbf{Case} $\mathbf{3.b.2}$: $N_{1} \gtrsim N$

We assume that $N_{4} \geq 1$. By (\ref{Eqn:EstMultRad2}) we have

\begin{equation}
\begin{array}{ll}
X & \lesssim \frac{N_{2}^{2(1-s)}}{N^{1-s} N_{1}^{1-s}} \|
\partial_{t} I u_{1} \|_{L_{t}^{6-} \Ivl L_{x}^{3+}}
\| I u_{2} \|_{L_{t}^{6} \Ivl L_{x}^{3}} \| I u_{3} \|_{L_{t}^{6} \Ivl L_{x}^{3}} \\
&  \| I u_{4} \|_{L_{t}^{2+} \Ivl L_{x}^{\infty -}} \\
& \lesssim N_{1}^{+} N_{4}^{-} \frac{N_{2}^{2(1-s)}}{N^{1-s}
N_{1}^{1-s}}  N_{1}^{\frac{1}{3}} N_{2}^{-\frac{2}{3}}
N_{3}^{-\frac{2}{3}}
\| \partial_{t} D^{- \left( \frac{1}{3}+ \right) } I u_{1} \|_{L_{t}^{6-} \Ivl L_{x}^{3+}} \\
& \| D^{1- \frac{1}{3} } I u_{2} \|_{L_{t}^{6} \Ivl L_{x}^{3}} \|
D^{1-\frac{1}{3}}  I u_{3} \|_{L_{t}^{6} \Ivl L_{x}^{3}}
\| D^{1-(1-)} I u_{4} \|_{L_{t}^{2+} \Ivl L_{x}^{\infty-}} \\
& \lesssim \frac{N_{2}^{-} N_{4}^{-}}{N^{1-}} Z^{4}(\tau)
\end{array}
\nonumber
\end{equation}
after directly creating $N_{1}^{+}$ and $N_{4}^{-}$. If $N_{4} < 1$
the proof is similar except that we indirectly create $N_{4}^{+}$.
We get (\ref{Eqn:AlmConsToProve}) after summation.

\end{itemize} 

\end{itemize} 

\item \textbf{Case} $\mathbf{4}$: $N_{4} \gtrsim N$

There are two subcases

\begin{itemize}

\item \textbf{Case} $\mathbf{4.a}$: $N_{1} \sim  N_{2}$

By (\ref{Eqn:EstMultRad2}) we have

\begin{equation}
\begin{array}{ll}
X & \lesssim \frac{N_{3}^{1-s}}{N^{1-s}} \frac{N_{4}^{1-s}}{N^{1-s}}
\| \partial_{t} I u_{1} \|_{L_{t}^{4} \Ivl L_{x}^{4}}
\| I u_{2} \|_{L_{t}^{4+} \Ivl L_{x}^{4-}} \| I u_{3} \|_{L_{t}^{4} \Ivl L_{x}^{4}} \\
& \| I u_{4} \|_{L_{t}^{4-} \Ivl L_{x}^{4+}} \\
& \lesssim  N_{2}^{-} N_{4}^{+} \frac{N_{3}^{1-s}}{N^{1-s}}
\frac{N_{4}^{1-s}} {N^{1-s}} N_{1}^{\frac{1}{2}}
\frac{1}{N_{2}^{\frac{1}{2}}} \frac{1}{N_{3}^{\frac{1}{2}}}
\frac{1}{N_{4}^{\frac{1}{2}}}
\| \partial_{t} D^{- \frac{1}{2} } I u_{1} \|_{L_{t}^{4} \Ivl L_{x}^{4}} \\
& \| D^{1- \left( \frac{1}{2} - \right) } I u_{2} \|_{L_{t}^{4+}
\Ivl L_{x}^{4-}}
\| D^{1-\frac{1}{2}}  I u_{3} \|_{L_{t}^{4} \Ivl L_{x}^{4}}  \| D^{1- \left( \frac{1}{2} + \right) } I u_{4} \|_{L_{t}^{4} \Ivl L_{x}^{4}} \\
& \lesssim \frac{N_{2}^{-}}{N^{1-}} Z^{4}(\tau)
\end{array}
\nonumber
\end{equation}
after directly creating $N_{2}^{-}$ and $N_{4}^{+}$.  We get
(\ref{Eqn:AlmConsToProve}) after summation.

\item \textbf{Case} $\mathbf{4.b}$: $N_{1} << N_{2}$

In this case we have $N_{2} \sim N_{3}$. There are two subcases

\begin{itemize}

\item \textbf{Case} $\mathbf{4.b.1}$: $N_{1} \gtrsim N$

By (\ref{Eqn:EstMultRad2}) we have

\begin{equation}
\begin{array}{ll}
X & \lesssim \frac{N_{2}^{2(1-s)}}{N^{2(1-s)}}
\frac{N_{4}^{1-s}}{N^{1-s}} \frac{N^{1-s}}{N_{1}^{1-s}}
 \| \partial_{t} I u_{1} \|_{L_{t}^{4} \Ivl L_{x}^{4}} \| I u_{2} \|_{L_{t}^{4} \Ivl L_{x}^{4}}
 \| I u_{3} \|_{L_{t}^{4} \Ivl L_{x}^{4}} \\
&  \| I u_{4} \|_{L_{t}^{4} \Ivl L_{x}^{4}} \\
& \lesssim \frac{N_{2}^{2(1-s)}}{N^{2(1-s)}}
\frac{N_{4}^{1-s}}{N^{1-s}} \frac{N^{1-s}}{N_{1}^{1-s}}
N_{1}^{\frac{1}{2}} \frac{1}{N_{2}^{\frac{1}{2}}}
\frac{1}{N_{3}^{\frac{1}{2}}} \frac{1}{N_{4}^{\frac{1}{2}}} \|
\partial_{t} D^{- \frac{1}{2} } I u_{1} \|_{L_{t}^{4} \Ivl
L_{x}^{4}}
\| D^{1- \frac{1}{2}  } I u_{2} \|_{L_{t}^{4} \Ivl L_{x}^{4}} \\
& \| D^{1-\frac{1}{2}}  I u_{3} \|_{L_{t}^{4} \Ivl L_{x}^{4}}  \| D^{1-  \frac{1}{2} } I u_{4} \|_{L_{t}^{4} \Ivl L_{x}^{4}} \\
& \lesssim \frac{N_{2}^{-}}{N^{1-}} Z^{4}(\tau)
\end{array}
\nonumber
\end{equation}
We get (\ref{Eqn:AlmConsToProve}) after summation.

\item \textbf{Case} $\mathbf{4.b.2}$: $N_{1} << N$

We assume that $N_{1} \geq 1$. We have

\begin{equation}
\begin{array}{ll}
X & \lesssim \frac{N_{2}^{2(1-s)}}{N^{2(1-s)}}
\frac{N_{4}^{1-s}}{N^{1-s}}
 \| \partial_{t} I u_{1} \|_{L_{t}^{4} \Ivl L_{x}^{4}} \| I u_{2} \|_{L_{t}^{4} \Ivl L_{x}^{4}}
 \| I u_{3} \|_{L_{t}^{4} \Ivl L_{x}^{4}} \\
&  \| I u_{4} \|_{L_{t}^{4} \Ivl L_{x}^{4}} \\
& \lesssim \frac{N_{2}^{2(1-s)}}{N^{2(1-s)}}
\frac{N_{4}^{1-s}}{N^{1-s}} N_{1}^{\frac{1}{2}}
\frac{1}{N_{2}^{\frac{1}{2}}} \frac{1}{N_{3}^{\frac{1}{2}}}
\frac{1}{N_{4}^{\frac{1}{2}}} \| \partial_{t} D^{- \frac{1}{2} } I
u_{1} \|_{L_{t}^{4} \Ivl L_{x}^{4}}
\| D^{1- \frac{1}{2}  } I u_{2} \|_{L_{t}^{4} \Ivl L_{x}^{4}} \\
& \| D^{1-\frac{1}{2}}  I u_{3} \|_{L_{t}^{4} \Ivl L_{x}^{4}}  \| D^{1-  \frac{1}{2} } I u_{4} \|_{L_{t}^{4} \Ivl L_{x}^{4}} \\
& \lesssim \frac{N_{1}^{-} N_{2}^{-} }{N^{1-}} Z^{4}(\tau)
\end{array}
\nonumber
\end{equation}
If $N_{1} <1$ the proof is similar except that we create $N_{1}^{+}$
instead of $N_{1}^{-}$. We get (\ref{Eqn:AlmConsToProve}) after
summation.

\end{itemize} 

\end{itemize} 

\end{itemize} 


\section{Proof of Almost Morawetz-Strauss inequality}
\label{sec:AlmMorIneq}

We prove proposition \ref{prop:AlmMorawetz} in this section. The
proof is divided into two steps

\begin{itemize}

\item \textbf{First Step}: Morawetz-Strauss inequality

We recall the proof of the Morawetz-Strauss inequality in \cite{mor,morstr}. We have the following identity

\begin{equation}
\begin{array}{l}
\left( \frac{x.\nabla u}{|x|} + \frac{u}{|x|} \right) \left( u_{tt}
- \triangle u + u^{3} \right)
= \partial_{t} \left( \frac{1}{|x|} \left( x.\nabla u + u  \right) \partial_{t}u \right) \\
+ div \left[ \frac{1}{|x|} \left(
\begin{array}{l}
-\frac{1}{2} \left( \partial_{t} u \right)^{2} -\left( x. \nabla u
\right) \nabla u +\frac{1}{2}|\nabla u|^{2}x -u \nabla u -
\frac{u^{2}}{2 |x|^{2}} x + \frac{1}{4} u^{4} x
\end{array}
\right)
\right] \\
+  \frac{1}{|x|} \left( |\nabla u|^{2} -\frac{\left( x. \nabla
u\right)^{2}}{|x|^{2}} \right) + \frac{u^{4}}{2 |x|}
\end{array}
\label{Eqn:Identity}
\end{equation}
and since u satisfies (\ref{Eqn:WaveEqRad}) we have after
integration

\begin{eqnarray}
2 \pi \int_{0}^{T} u^{2}(t,0) dt + \int_{0}^{T}
\int_{\mathbb{R}^{3}} \frac{u^{4}(t,x)}{2 |x|} dx \, dt
& = & - \int_{\mathbb{R}^{3}} \left( \frac{\nabla u (T,x).x}{|x|} + \frac{u(T,x)}{|x|} \right) \partial_{t} u (T,x) dx \nonumber \\
& & + \int_{\mathbb{R}^{3}} \left( \frac{\nabla u (0,x).x}{|x|} + \frac{u(0,x)}{|x|} \right) \partial_{t} u (0,x) dx \nonumber \\
& &
\nonumber
\end{eqnarray}
Now we apply the basic inequality $|ab| \leq \frac{|a|^{2}}{2} +
\frac{|b|^{2}}{2}$ to the right hand side of the integral and we get

\begin{eqnarray}
\int_{0}^{T} \int_{\mathbb{R}^{3}} \frac{u^{4}(t,x)}{2 |x|} dx \, dt
& \leq & \frac{1}{2}  \int_{\mathbb{R}^{3}} \left( \frac{\nabla u
(T,x).x}{|x|} + \frac{u(T,x)}{|x|} \right)^{2} +
(\partial_{t} u)^{2}(T,x) dx  \nonumber \\
& & + \frac{1}{2}  \int_{\mathbb{R}^{3}} \left( \frac{\nabla u
(0,x).x}{|x|} + \frac{u(0,x)}{|x|} \right)^{2} +
(\partial_{t} u)^{2}(0,x) dx  \nonumber \\
& & \label{Eqn:MorCauch}
\end{eqnarray}
We also notice that

\begin{eqnarray}
 \left( \frac{\nabla u.x}{|x|} + \frac{u}{|x|} \right)^{2} & = & \frac{ \left( \nabla u.x \right)^{2}}{|x|^{2}}
 + div \left( \frac{u^{2}}{|x|^{2}} x \right) \nonumber \\
 & \leq & | \nabla u |^{2} +  div \left( \frac{u^{2}}{|x|^{2}} x \right)
\label{Eqn:PtwiseEquality}
\end{eqnarray}
We plug (\ref{Eqn:PtwiseEquality}) into (\ref{Eqn:MorCauch}). We get
the Morawetz-Strauss's inequality

\begin{eqnarray}
\int_{0}^{T} \int_{ \mathbb{R}^{3}} \frac{u^{4}(t,x)}{ |x|} dx \, dt
& \leq & 2 \left( E(u(T))+ E(u(0)) \right)
\end{eqnarray}

\item \textbf{Second Step}: Almost Morawetz-Strauss's inequality.
We substitute $u$ for $Iu$ in (\ref{Eqn:Identity}) and we proceed
similarly. We get

\begin{equation}
\begin{array}{ll}
\int_{0}^{T} \int_{\mathbb{R}^{3}} \frac{|Iu |^{4}(t,x)}{ |x|} dx \,
dt -2 \left( E \left( Iu(T) \right) + E \left( Iu(0)  \right)
\right)
& \leq |R_{1}(T) + R_{2}(T)| \\
& \leq |R_{1}(T)| + |R_{2}(T)|
\end{array}
\label{Eqn:IntermMor}
\end{equation}

\end{itemize}

\section{Proof of the integral estimates}
\label{subsec:RemainderMor}

We are interested in proving proposition \ref{prop:RemainderMor} in
this section. In what follows we also assume that $J=[0, \, \tau]$:
the reader can check after reading the proof that the other cases
can be reduced to that one.

Plancherel formula yields

\begin{equation}
\begin{array}{ll}
R_{1}(\tau) & = \int_{0}^{\tau} \int_{\xi_{1}+... + \xi_{4}=0}
\mu(\xi_{2}, \xi_{3}, \xi_{4}) \widehat{\frac{\nabla I u.x}{|x|}}
(t,\xi_{1}) \widehat{I u} (t,\xi_{2}) \widehat{I u}(t,\xi_{3})
\widehat{I u} (t,\xi_{4}) d \xi_{2}...d \xi_{4} \, dt
\end{array}
\nonumber
\end{equation}
and

\begin{equation}
\begin{array}{ll}
R_{2}(\tau) & = \int_{0}^{\tau} \int_{\xi_{1}+...+\xi_{4}}
\mu(\xi_{2},\xi_{3},\xi_{4}) \widehat{\frac{Iu}{|x|}} (t,\xi_{1})
\widehat{I u} (t,\xi_{2}) \widehat{I u}(t,\xi_{3}) \widehat{I u}
(t,\xi_{4}) d \xi_{2}...d \xi_{4} \, dt
\end{array}
\nonumber
\end{equation}
with $\mu$ defined in (\ref{Eqn:DfnMu}). It suffices to prove

\begin{equation}
\begin{array}{ll}
\left| \int_{0}^{\tau} \int_{\xi_{1}+... + \xi_{4}=0}
\mu(\xi_{2},\xi_{3},\xi_{4}) \widehat{\frac{\nabla I u.x}{|x|}}
(t,\xi_{1}) \widehat{I u} (t,\xi_{2}) \widehat{I u}(t,\xi_{3})
\widehat{I u} (t,\xi_{4}) d \xi_{2}...d \xi_{4} \, dt \right|   &
\lesssim \frac{Z^{4}(\tau)}{N^{1-}}
\end{array}
\label{Eqn:R1Four}
\end{equation}
and

\begin{equation}
\begin{array}{ll}
\left| \int_{0}^{\tau} \int_{\xi_{1}+... + \xi_{4}=0}
\mu(\xi_{2},\xi_{3},\xi_{4}) \widehat{\frac{I u}{|x|}} (t,\xi_{1})
\widehat{I u} (t,\xi_{2}) \widehat{I u}(t,\xi_{3}) \widehat{I u}
(t,\xi_{4}) d \xi_{2}...d \xi_{4} \, dt \right|   & \lesssim
\frac{Z^{4}(\tau)}{N^{1-}}
\end{array}
\label{Eqn:R2Four}
\end{equation}
We perform a Paley-Littlewood decomposition to prove
(\ref{Eqn:R1Four}) and (\ref{Eqn:R2Four}). Let $u_{i}:=P_{N_{i}} u$,
$i \in \{2,\, ..., \,4\}$, $\left( \frac{\nabla I u \cdot x}{|x|}
\right)_{1}:= P_{N_{1}} \left( \frac{\nabla I u \cdot x}{|x|}
\right) $ and $ \left( \frac{I u}{|x|} \right)_{1} := P_{N_{1}}
\left( \frac{Iu}{|x|} \right)$.

\begin{equation}
\begin{array}{ll}
X_{1} & = \left| \int_{0}^{\tau} \int_{\xi_{1}+... + \xi_{4}=0}
\mu(\xi_{2},\xi_{3},\xi_{4}) \widehat{\left(\frac{\nabla I u.x}{|x|}
\right)_{1}} (t,\xi_{1}) \widehat{I u_{2}} (t,\xi_{2}) \widehat{I
u_{3}}(t,\xi_{3}) \widehat{I u_{4}} (t,\xi_{4}) d \xi_{2}...d
\xi_{4} \, dt \right|
\end{array}
\nonumber
\end{equation}
and

\begin{equation}
\begin{array}{ll}
X_{2} & = \left| \int_{0}^{\tau} \int_{\xi_{1}+... + \xi_{4}=0} \mu
(\xi_{2}, \xi_{3}, \xi_{4}) \widehat{\left( \frac{I u}{|x|}
\right)_{1}} (t,\xi_{1}) \widehat{I u_{2}} (t,\xi_{2}) \widehat{I
u_{3}}(t,\xi_{3}) \widehat{I u_{4}} (t,\xi_{4}) d \xi_{2}...d
\xi_{4} \, dt \right|
\end{array}
\nonumber
\end{equation}
Notice that by Bernstein inequality, H\"older inequality, Plancherel
theorem and (\ref{Eqn:HardyIneq}) we have

\begin{equation}
\begin{array}{ll}
\left\| \left( \frac{\nabla I u \cdot x}{|x|} \right)_{1} \right\|_{
L_{t}^{\infty -} L_{t}^{2+}} & \lesssim N_{1}^{+}  \left\|
\frac{\nabla I u \cdot
x}{|x|} \right\|_{L_{t}^{\infty} L_{x}^{2}} \\
& \lesssim N_{1}^{+} \| \nabla I u \|_{L_{t}^{\infty} L_{x}^{2}} \\
& \lesssim N_{1}^{+} \| D I u \|_{L_{t}^{\infty} L_{x}^{2}}
\end{array}
\label{Eqn:1Mor}
\end{equation}
and

\begin{equation}
\begin{array}{ll}
\left\| \left( \frac{ I u}{|x|} \right)_{1}
\right\|_{L_{t}^{\infty-} L_{x}^{2+}} & \lesssim N_{1}^{+} \left\|
\frac{Iu}{|x|}
\right\|_{L_{t}^{\infty} L_{x}^{2}} \\
& \lesssim N_{1}^{+} \| D I u \|_{L_{t}^{\infty} L_{x}^{2}}
\end{array}
\label{Eqn:2Mor}
\end{equation}

If $p_{j} \in [1, \, \infty]$ and $q_{j} \in (1, \, \infty)$, $j \in
\{2, \, ..., \, 4\}$ such that $ \frac{1}{(\infty-)} +
\sum_{j=2}^{4} \frac{1}{p_{j}} = 1$, $ \frac{1}{2+} + \sum_{j=2}^{4}
\frac{1}{q_{j}} =1$, $(p_{j},q_{j})$ -$m_{j}$ wave admissible for
some $m_{j}^{'} \, s$ such that $0 \leq m_{j} < 1$ and
$\frac{1}{p_{j}} + \frac{1}{q_{j}}=\frac{1}{2}$ then we have by the
methodology explained in the proof of proposition
\ref{prop:EstNrjRad}

\begin{equation}
\begin{array}{ll}
X_{1}  & \lesssim B \left( N_{2},N_{3},N_{4} \right)
 \left\| \left( \frac{\nabla Iu \cdot x}{|x|} \right)_{1} \right\|_{L_{t}^{\infty-} \Ivl L_{x}^{2+}}
\| I u_{2} \|_{L_{t}^{p_{2}} \Ivl L_{x}^{q_{2}}} ... \| I u_{4}
\|_{L_{t}^{p_{4}} \Ivl L_{x}^{q_{4}}}
\end{array}
\nonumber
\end{equation}
and

\begin{equation}
\begin{array}{ll}
X_{2} & \lesssim B \left( N_{2},N_{3},N_{4} \right) \left\| \left(
\frac{I u}{|x|} \right)_{1} \right\|_{L_{t}^{\infty -} \Ivl
L_{x}^{2+}} \| I u_{2} \|_{L_{t}^{p_{2}} \Ivl L_{x}^{q_{2}}} ... \|
I u_{4} \|_{L_{t}^{p_{4}} \Ivl L_{x}^{q_{4}}}
\end{array}
\nonumber
\end{equation}
By symmetry we can assume that $N_{2} \geq N_{3} \geq N_{4}$. There
are different cases

\begin{itemize}

\item \textbf{Case} $\mathbf{1}$: $N >> N_{2} \geq N_{3}$. In this case
$X_{1}=0$ and $X_{2}=0$ since $\mu=0$.

\item \textbf{Case} $\mathbf{2}$: $N_{2} \gtrsim N >> N_{3}$
By (\ref{Eqn:IndirectRes}), (\ref{Eqn:EstMultRad1}),
(\ref{Eqn:1Mor}) and (\ref{Eqn:2Mor}) we have

\begin{equation}
\begin{array}{ll}
X_{1} & \lesssim \frac{N_{3}}{N_{2}} \left\| \left( \frac{\nabla
Iu.x}{|x|} \right)_{1} \right\|_{L_{t}^{\infty-} \Ivl  L_{x}^{2+}}
\| I u_{2} \|_{L_{t}^{\infty} \Ivl L_{x}^{2}} \| I u_{3}
\|_{L_{t}^{2+} \Ivl
L_{x}^{\infty-}} \| I u_{4} \|_{L_{t}^{2+} \Ivl L_{x}^{\infty -}} \\
& \lesssim \frac{N_{3}}{N_{2}} N_{1}^{+} \frac{1}{N_{2}} N_{3}^{+}
N_{4}^{+} \| D I u \|_{L_{t}^{\infty} \Ivl L_{x}^{2}} \| D I u_{2}
\|_{L_{t}^{\infty} \Ivl L_{x}^{2}} \| D^{1-(1-)} I u_{3} \|_{L_{t}^{2+} \Ivl
L_{x}^{\infty}} \\
& \| D^{1-(1-)} I u_{4} \|_{L_{t}^{2+} \Ivl
L_{x}^{\infty-}} \\
& \lesssim \frac{ N_{2}^{--} N_{4}^{+}}{N^{1-}} Z^{4}(\tau)
\end{array}
\nonumber
\end{equation}
and

\begin{equation}
\begin{array}{ll}
X_{2} & \lesssim \frac{N_{3}}{N_{2}} \left\| \left( \frac{Iu}{|x|}
\right)_{1} \right\|_{L_{t}^{\infty-} \Ivl  L_{x}^{2+}} \| I u_{2}
\|_{L_{t}^{\infty} \Ivl L_{x}^{2}} \| I u_{3} \|_{L_{t}^{2+} \Ivl
L_{x}^{\infty-}} \| I u_{4} \|_{L_{t}^{2+} \Ivl L_{x}^{\infty -}} \\
& \lesssim \frac{N_{3}}{N_{2}} N_{1}^{+} \frac{1}{N_{2}} N_{3}^{+}
N_{4}^{+} \| D I u \|_{L_{t}^{\infty} \Ivl L_{x}^{2}} \| D I u_{2}
\|_{L_{t}^{\infty} \Ivl L_{x}^{2}} \| D^{1-(1-)} I u_{3} \|_{L_{t}^{2+} \Ivl
L_{x}^{\infty}} \\
& \| D^{1-(1-)} I u_{4} \|_{L_{t}^{2+} \Ivl
L_{x}^{\infty-}} \\
& \lesssim \frac{ N_{2}^{--} N_{4}^{+}}{N^{1-}} Z^{4}(\tau)
\end{array}
\nonumber
\end{equation}

\item \textbf{Case} $\mathbf{3}$: $N_{3} \gtrsim N >> N_{4}$

There are two subcases:

\begin{itemize}

\item \textbf{Case} $\mathbf{3.a}$: $N_{1} \sim N_{2}$

By (\ref{Eqn:IndirectRes}), (\ref{Eqn:EstMultRad2}) and
(\ref{Eqn:1Mor})

\begin{equation}
\begin{array}{ll}
X_{1} & \lesssim \frac{N_{3}^{1-s}}{N^{1-s}} \left\| \left(
\frac{\nabla Iu.x}{|x|} \right)_{1} \right\|_{L_{t}^{\infty-} \Ivl
L_{x}^{2+}} \| I u_{2} \|_{L_{t}^{\infty} \Ivl L_{x}^{2}} \| I u_{3}
\|_{L_{t}^{2+} \Ivl
L_{x}^{\infty-}} \| I u_{4} \|_{L_{t}^{2+} \Ivl L_{x}^{\infty -}} \\
& \lesssim \frac{N_{3}^{1-s}}{N^{1-s}} N_{1}^{+} \frac{1}{N_{2}}
N_{3}^{+} N_{4}^{+} \| D I u \|_{L_{t}^{\infty} Ivl L_{x}^{2}} \| D I
u_{2} \|_{L_{t}^{\infty} \Ivl L_{x}^{2}} \| D^{1-(1-)} I u_{3}
\|_{L_{t}^{2+} \Ivl L_{x}^{\infty}} \\
& \| D^{1-(1-)} I u_{4} \|_{L_{t}^{2+} \Ivl
L_{x}^{\infty-}} \\
& \lesssim \frac{ N_{2}^{--} N_{4}^{+}}{N^{1-}} Z^{4}(\tau)
\end{array}
\label{Eqn:Case3a}
\end{equation}
Similarly we get $X_{2} \lesssim \frac{ N_{2}^{--}
N_{4}^{+}}{N^{1-}} Z^{4}(\tau)$ after substituting $X_{1}$, $
\left\| \left( \frac{\nabla Iu.x}{|x|} \right)_{1}
\right\|_{L_{t}^{\infty-} \Ivl L_{x}^{2+}} $ for $X_{2}$, $ \left\|
\left( \frac{Iu}{|x|} \right)_{1} \right\|_{L_{t}^{\infty-} \Ivl
L_{x}^{2+}} $ respectively in (\ref{Eqn:Case3a}).

\item \textbf{Case} $\mathbf{3.b}$: $N_{1} << N_{2}$

There are two subcases

\begin{itemize}

\item \textbf{Case} $\mathbf{3.b.1}$ $N_{1} << N$

\begin{equation}
\begin{array}{ll}
X_{1} & \lesssim \frac{N_{2}^{2(1-s)}}{N^{2(1-s)}} \left\| \left(
\frac{\nabla Iu.x}{|x|} \right)_{1} \right\|_{L_{t}^{\infty-} \Ivl
L_{x}^{2+}} \| I u_{2} \|_{L_{t}^{\infty} \Ivl L_{x}^{2}} \| I u_{3}
\|_{L_{t}^{2+} \Ivl
L_{x}^{\infty-}} \| I u_{4} \|_{L_{t}^{2+} \Ivl L_{x}^{\infty -}} \\
& \lesssim \frac{N_{2}^{2(1-s)}}{N^{2(1-s)}} N_{1}^{+}
\frac{1}{N_{2}} N_{3}^{+} N_{4}^{+} \| D I u \|_{L_{t}^{\infty} \Ivl
L_{x}^{2}} \| D I u_{2} \|_{L_{t}^{\infty} \Ivl L_{x}^{2}} \| D^{1-(1-)}
I u_{3} \|_{L_{t}^{2+} \Ivl L_{x}^{\infty}} \\
& \| D^{1-(1-)} I u_{4} \|_{L_{t}^{2+} \Ivl
L_{x}^{\infty-}} \\
& \lesssim \frac{ N_{1}^{+} N_{2}^{---} N_{4}^{+}}{N^{1-}}
Z^{4}(\tau)
\end{array}
\nonumber
\end{equation}
Similarly $ X_{2} \lesssim \frac{ N_{1}^{+} N_{2}^{---}
N_{4}^{+}}{N^{1-}} Z^{4}(\tau)$.

\item \textbf{Case} $\mathbf{3.b.2}$ $N_{1} \gtrsim N$

\begin{equation}
\begin{array}{ll}
X_{1} & \lesssim \frac{N_{2}^{2(1-s)}}{N^{1-s} N_{1}^{1-s}} \left\|
\left( \frac{\nabla Iu.x}{|x|} \right)_{1} \right\|_{L_{t}^{\infty-}
\Ivl L_{x}^{2+}} \| I u_{2} \|_{L_{t}^{\infty} \Ivl L_{x}^{2}} \| I
u_{3} \|_{L_{t}^{2+} \Ivl
L_{x}^{\infty-}} \| I u_{4} \|_{L_{t}^{2+} \Ivl L_{x}^{\infty -}} \\
& \lesssim \frac{N_{2}^{2(1-s)}}{N^{1-s} N_{1}^{1-s}} N_{1}^{+}
\frac{1}{N_{2}} N_{3}^{+} N_{4}^{+} \| D I u \|_{L_{t}^{\infty} \Ivl
L_{x}^{2}} \| D I u_{2} \|_{L_{t}^{\infty} \Ivl L_{x}^{2}} \| D^{1-(1-)}
I u_{3} \|_{L_{t}^{2+} \Ivl L_{x}^{\infty}} \\
& \| D^{1-(1-)} I u_{4}
\|_{L_{t}^{2+} \Ivl
L_{x}^{\infty-}} \\
& \lesssim \frac{ N_{2}^{--} N_{4}^{+}}{N^{1-}} Z^{4}(\tau)
\end{array}
\nonumber
\end{equation}
Similarly $X_{2} \lesssim \frac{ N_{2}^{--} N_{4}^{+}}{N^{1-}}
Z^{4}(\tau)$.

\end{itemize}

\end{itemize}

\item \textbf{Case} $\mathbf{4}$: $N_{4} \gtrsim N$

There are two subcases

\begin{itemize}

\item \textbf{Case} $\mathbf{4.a}$: $N_{1} \sim N_{2}$

\begin{equation}
\begin{array}{ll}
X_{1} & \lesssim \frac{N_{3}^{1-s}}{N^{1-s}}
\frac{N_{4}^{1-s}}{N^{1-s}} \left\| \left( \frac{ \nabla Iu.x}{|x|}
\right)_{1} \right\|_{L_{t}^{\infty -} \Ivl L_{x}^{2+}} \left\| I u_{2}
\right \|_{L_{t}^{\infty}([0,\tau]) L_{x}^{2}} \| I u_{3}
\|_{L_{t}^{2+}([0, \tau]) L_{x}^{\infty-}} \|
I u_{4} \|_{L_{t}^{2+} \Ivl L_{x}^{\infty-}} \\
& \lesssim \frac{N_{3}^{1-s}}{N^{1-s}} \frac{N_{4}^{1-s}}{N^{1-s}}
N_{1}^{+} \frac{1}{N_{2}} N_{3}^{+} N_{4}^{+} \| D I u
\|_{L_{t}^{\infty} \Ivl L_{x}^{2}} \| D I u_{2} \|_{L_{t}^{\infty} \Ivl
L_{x}^{2}} \| D^{1-(1-)} I u_{3} \|_{L_{t}^{2+} \Ivl L_{x}^{\infty}} \\
& \| D^{1-(1-)} I u_{4} \|_{L_{t}^{2+} \Ivl L_{x}^{\infty-}} \\
& \lesssim \frac{N_{2}^{-}} {N^{1-}} Z^{4}(\tau)
\end{array}
\nonumber
\end{equation}
Similarly $X_{2} \lesssim \frac{N_{2}^{-}}{N^{1-}} Z^{4}(\tau)$.

\item \textbf{Case} $\mathbf{4.b}$: $N_{1}<< N_{2}$. There are two subcases

\begin{itemize}

\item \textbf{Case} $\mathbf{4.b.1}$: $N_{1} \gtrsim N$. We have

\begin{equation}
\begin{array}{ll}
X_{1} & \lesssim \frac{N_{2}^{2(1-s)}}{N^{2(1-s)}}
\frac{N_{4}^{1-s}}{N^{1-s}} \frac{N^{1-s}}{N_{1}^{1-s}} \left\|
\left( \frac{ \nabla Iu.x}{|x|} \right)_{1} \right\|_{L_{t}^{\infty
-} \Ivl L_{x}^{2+}} \left\| I u_{2} \right \|_{L_{t}^{\infty}([0,\tau])
L_{x}^{2}} \| I u_{3} \|_{L_{t}^{2+}([0, \tau]) L_{x}^{\infty-}} \|
I u_{4} \|_{L_{t}^{2+} \Ivl L_{x}^{\infty-}} \\
& \lesssim \frac{N_{2}^{2(1-s)}}{N^{2(1-s)}}
\frac{N_{4}^{1-s}}{N^{1-s}} \frac{N^{1-s}}{N_{1}^{1-s}}  N_{1}^{+}
\frac{1}{N_{2}} N_{3}^{+} N_{4}^{+} \| D I u \|_{L_{t}^{\infty} \Ivl
L_{x}^{2}} \| D I u_{2} \|_{L_{t}^{\infty} \Ivl L_{x}^{2}} \| D^{1-(1-)}
I u_{3} \|_{L_{t}^{2+} \Ivl L_{x}^{\infty}} \\
& \| D^{1-(1-)} I u_{4} \|_{L_{t}^{2+} \Ivl L_{x}^{\infty-}} \\
& \lesssim \frac{N_{2}^{-}}{N^{1-}} Z^{4}(\tau)
\end{array}
\nonumber
\end{equation}
Similarly $X_{2} \lesssim \frac{N_{2}^{-}}{N^{1-}} Z^{4}(\tau)$.

\item \textbf{Case} $\mathbf{4.b.2}$: $N_{1} << N$. We have

\begin{equation}
\begin{array}{ll}
X_{1} & \lesssim \frac{N_{2}^{2(1-s)}}{N^{2(1-s)}}
\frac{N_{4}^{1-s}}{N^{1-s}} \left\| \left( \frac{\nabla Iu.x}{|x|}
\right)_{1} \right\|_{L_{t}^{\infty -} \Ivl L_{x}^{2+}} \left\| I u_{2}
\right \|_{L_{t}^{\infty}([0,\tau]) L_{x}^{2}} \| I u_{3}
\|_{L_{t}^{2+}([0, \tau]) L_{x}^{\infty-}} \|
I u_{4} \|_{L_{t}^{2+} \Ivl L_{x}^{\infty-}} \\
& \lesssim \frac{N_{2}^{2(1-s)}}{N^{2(1-s)}}
\frac{N_{4}^{1-s}}{N^{1-s}} N_{1}^{+} \frac{1}{N_{2}} N_{3}^{+}
N_{4}^{+} \| D I u \|_{L_{t}^{\infty} \Ivl L_{x}^{2}} \| D I u_{2}
\|_{L_{t}^{\infty} \Ivl L_{x}^{2}} \| D^{1-(1-)} I u_{3} \|_{L_{t}^{2+} \Ivl
L_{x}^{\infty}} \\
& \| D^{1-(1-)} I u_{4} \|_{L_{t}^{2+} \Ivl L_{x}^{\infty-}} \\
& \lesssim \frac{N_{1}^{+} N_{2}^{--}}{N^{1-}} Z^{4}(\tau)
\end{array}
\end{equation}
Similarly $X_{2} \lesssim \frac{N_{1}^{+} N_{2}^{--}}{N^{1-}}
Z^{4}(\tau)$.

\end{itemize}
\end{itemize}
\end{itemize}
We get (\ref{Eqn:R1Four}) and (\ref{Eqn:R2Four}) after summation.


\begin{thebibliography}{99}


\bibitem{bahchemin} H.Bahouri and Jean-Yves Chemin; \emph{On global well-posedness for defocusing cubic wave equation}, to appear
Internat. Math. Res. Notices

\bibitem{bourg} J. Bourgain; \emph{Refinement of Strichartz inequality
and applications to $2D-NLS$ with critical nonlinearity}, Internat.
Math. Res. Notices, 5 (1998), 253-283

\bibitem{caz} T. Cazenave; \emph{Semilinear Schr\"odinger Equations}, Amer. Math. Soc., Providence, Rhode Island, 2003

\bibitem{coifmey2} R. R Coifman and Y. Meyer; \emph{Commutateurs d'integrales singuli\`eres et op\'erateurs multilin\'eaires
}, Ann. Inst. Fourier (Grenoble) 28, pp. 177-202 [1978]


\bibitem{almckstt} J.Colliander, M.Keel, G.Staffilani, H.Takaoka, T.Tao;
\emph{Almost conservation laws and global rough solutions to a
nonlinear Schr\"odinger equation}, Math. Res. Letters 9 (2002), pp.
659-682

\bibitem{morckstt} J.Colliander, M.Keel, G.Staffilani, H.Takaoka, T.Tao;
\emph{Global well-posedness and scattering for rough solutions of a
nonlinear Schr\"odinger equation on $\mathbb{R}^{3}$ }, CPAM 57
(2004), pp 1-34


\bibitem{gallagplanch} I. Gallagher and F. Planchon; \emph{On global solutions to a dofocusing
semi-linear wave equation}, Revista Mathematic$\acute{a}$
Iberoamericana, 19, 2003, pp. 161-177

\bibitem{ginebvelo} J. Ginebre and G. Velo; \emph{Generalized Strichartz inequalities for the
wave equation}, Jour. Func. Anal. 133 (1995), pp 50-68

\bibitem{kenponcevega} C. E. Kenig, G. Ponce and L. Vega; \emph{Global well-posedness for
semi-linear wave equations}, Communications in partial differential
equations 25 (2000), pp. 1741-1752

\bibitem{klaintat} S. Klainerman and D. Tataru; \emph{On the optimal local regularity for
the Yang-Mills equations in $\mathbb{R}^{4+1}$}, Journal of the
American Mathematical Society 12 (1999), pp. 93-116

\bibitem{lindsogge} H. Lindbald and C. Sogge; \emph{On existence and scattering with
minimal regularity for semilinear wave equations}; J. Funct. Anal
130 (1995), pp. 357-426

\bibitem{mor} C. Morawetz; \emph{Time decay for the nonlinear Klein-Gordon equation}; Proc. Roy. Soc. A 306, 1968, pp.291-296

\bibitem{morstr} C. Morawetz and W. Strauss; \emph{Decay and scattering of solutions of a
nonlinear relativistic wave equation}; Comm. Pure Appl. Math. 25
(1972), pp 1-31




\end{thebibliography}
\end{document}